\documentclass[10pt]{article}
 \usepackage{latexsym,amsmath,amssymb,amsbsy,amstext,amscd,amsfonts}
 \usepackage{graphics,graphicx}
 \usepackage{color}
 \usepackage{tabularx}
 \usepackage{multirow}
\usepackage{booktabs}

\date{\today}

\title{Efficient pseudo-spectral solvers for the PKN model of hydrofracturing.}

\author{Michal Wrobel$^{(2,1)}$, Gennady Mishuris$^{(1,2)}$,
\\
{\it $^{(1)}$ Institute of Mathematical and Physical Sciences,
Aberystwyth University, }
\\ {\it Ceredigion SY23 3BZ, Wales U.K.,}
\\{\it $^{(2)}$ Eurotech Sp. z o.o.,  }
\\ {\it ul. Wojska Polskiego 3, 39-300 Mielec, Poland}}

\begin{document}

\maketitle

\begin{abstract}
In the paper,  a novel algorithm employing
pseudo-spectral approach is developed for the PKN mo\-del of hydrofracturing.
The respective solvers compute both the solution and its temporal derivative. In comparison with conventional solvers, they
demonstrate significant cost effectiveness in terms of balance between the accuracy of computations and densities of the temporal and spatial meshes.
Various fluid flow regimes are considered.
\end{abstract}

\section{Introduction and preliminary results}
Hydraulic fracturing is a widely used method for stimulation of hydrocarbons reservoirs. This technology has been known and
successfully applied for a few decades \cite{Khristianovic,Hubbert,Crittendon}. Recently it has been
revived, due to economical reasons, as a basic technique for
exploitation of unconventional deposits of oil and gas. The
phenomenon of a fluid driven fracture propagating in a brittle
medium is also present in many natural processes (e.g. magma
driven dykes -- \cite{Rubin}, subglacial drainage of water -- \cite{Rice}).

Throughout the years, starting from the pioneering works of Sneddon and Elliot \cite{Sneddon}, Khristianovic and Zheltov \cite{Khristianovic}, Perkins and Kern \cite{Perkins}, Geertsma and de Klerk \cite{Geertsma}, and Nordgren \cite{Nordgren}, various models of hydrofracturing have been
formulated and used in applications. A broad review of the topic can
be found in \cite{Kovalyshen,Kovalyshen_1,Linkov_3,Kusmierczyk}. Together with increasing complexity of the
models describing this multhiphysics process, the computational
techniques have been continuously enhanced. A comprehensive survey on the
algorithms and numerical methods used in hydrofracturing simulation
can be found in \cite{Adachi_2007,Economides2000}.

Responding to the recent demand, an increasing stream of publications
have appeared concerning additional information on seismic
events, shear stresses in the rock formation, multifracturing and
others and their implementation into the  solvers \cite{Zhang,Moos,Olson,Dobroskok}. Also, a
considerable effort has been made to improve the existing algorithms by
incorporating new efficient numerical techniques \cite{Lecampion_Det,Linkov_3,Pierc_Det,Pierc_Det_2}.

The main computational challenges associated with the modelling of
hydraulic fractures are: a) strong nonlinearity resulting from the
coupling between the solid and fluid phases, b) singularity of the
gradients of the physical fields near the crack tip, c) moving boundaries, d) degeneration of
the governing equations at the crack tip, multiscaling and others.

To achieve the maximal possible
efficiency of numerical simulations, the computational algorithms should be formulated in {\it
proper} variables accounting for all the problem
peculiarities \cite{Linkov_3}. As a result,  they allow one to reduce the
volume of processed data, which is especially important when dealing
with complex geometries and/or multifracturing.

The analysis presented in this paper is devoted to the PKN model of
hydrofracturing. This model contains all the peculiarities mentioned
above, except for the non-local relation for
the fluid-solid coupling. Although we restrict our interest only to a single fracture, the developed algorithms, thanks to their robustness, can be
successfully applied to model a system of cracks.

The numerical analysis of the problem should be backdated to Nordgren
 \cite{Nordgren} who extended the Perkins and Kern model
\cite{Perkins} to account for the fluid loss effect and fracture
volume change. As a result, the crack length was determined  as
part of the solution. The author proposed a finite difference scheme to
solve the problem, which is in fact equivalent to the finite volume
(FV) method.

Further development of the PKN formulation was done by Kemp
\cite{Kemp}, who (a) implemented the
specific boundary condition at the moving crack tip into the FV scheme,  (b)
incorporated asymptotic behaviour of the solution near the crack tip
in a  special tip element, (c) indirectly used the fourth power of the crack opening ($w^4$) as a new dependent
variable, instead of the
crack opening itself. For the early-time asymptotic model Kemp
proposed a power series solution, presenting its first four terms.

The recent paper by Kovalyshen and Detournay\cite{Kovalyshen} has extended most of Kemp's
results, incorporating all
information on the PKN model available to date. They present various asymptotics,
complete analytical solution for an impermeable rock (directly extending
the results from \cite{Kemp} from four leading terms to an infinite
series representation), FV algorithm with a special tip element and a numerical benchmark for the Carter leak-off,
linking the results to the scaling approach developed in \cite{Adachi_2007b,Detournay,Mitchell_2007,Mitchell_2007b}.

In \cite{Linkov_1}-\cite{Linkov_3}, the PKN model was reformulated by Linkov to improve the efficiency and stability of computations
by (i) introducing {\it proper} dependent variables (cubed fracture opening, $w^3$), (ii)
utilizing the speed equation and (iii)
by imposing a modified boundary condition at a
small distance behind the crack tip ($\varepsilon$-regularisation).
Additionally, the analytical
solution for an impermeable rock was evaluated for the new dependent variable in a form of rapidly converging series in \cite{Linkov_3}. Moreover, the
author highlighted in \cite{Linkov_1} that
numerical schemes exploiting a fixed position of the crack tip during the iterations may become ill-posed.

In \cite{M_W_L} and \cite{Kusmierczyk} the
$\varepsilon$-regularisation technique was further enhanced by (i) appropriate adaptation of the speed equation
to the chosen numerical scheme and (ii) improved way of imposing of the regularized boundary condition. A detailed
discussion on various aspects of application of implicit and
explicit numerical schemes was provided.

In this paper we are presenting a novel algorithm based on the
pseudo-spectral approach. Namely, we propose an efficient numerical algorithm to solve a specific self-similar
problem and extend the results
to the general (transient) formulation. Since, the
integration schemes used in the algorithm incorporate the exact
boundary conditions at the crack tip, no regularization technique is
necessary. The most accurate two points representation of the temporal derivative is used to guarantee an optimal algorithm
performance. Finally, two solvers are developed which show their robustness and
stability. They both demonstrate high cost effectiveness in
terms of the relationship between the volume of the processed data and
the accuracy of computations. Moreover, additionally to the crack opening and length,
the temporal derivative of the former and the crack tip velocity are automatically returned as components of the problem solution.

\subsection{Problem formulation}

Let us consider a symmetrical crack of length $2l$ situated in
the plane $x\in [-l,l]$. The crack is fully filled by a
Newtonian liquid injected at the middle point ($x=0$) with a known rate $q_0(t)$. Note here, that
the crack length evolution, $l=l(t)$, is the result of fluid flow inside the fracture. Due to
the symmetry of the problem, one can restrict the analysis to the half of the crack $x\in [0,l(t)]$.

The classic mathematical formulation of the PKN model of hydrofracturing
was given in \cite{Nordgren}. Below we present a system of equations
constituting the model. The mass conservation principle is expressed
by the continuity equation:
\begin{equation}
\label{continuos_1} \frac{\partial w}{\partial t}+\frac{\partial
q}{\partial x}+q_l=0,\quad t \geq t_0,\quad 0\leq x \leq l(t),
\end{equation}
while the Poiseuille equation describes the flow in a narrow
channel. In the case of a Newtonian fluid, it is written in the following
form:
\begin{equation}\label{Poiseulle_1}
q=-\frac{1}{M}w^3\frac{\partial p}{\partial x}.
\end{equation}
Here  $w=w(t,x)$ stands for the crack opening, $q=q(t,x)$ is the
fluid flow rate,   $p=p(t,x)$ ($p=p_{f}-\sigma_0$, $\sigma_0$ -
confining stress) refers to the net fluid pressure. The constant $M$,
involved in the Poiseuille equation, is computed as $M=12{\mu}$,
where $\mu$ denotes the dynamic viscosity (see for example
\cite{Adachi_2002}). The function $q_l=q_l(t,x)$ from \eqref{continuos_1} is the volumetric rate of fluid loss to
formation in the direction perpendicular to the crack surfaces per unit
length of the fracture. This function is usually assumed to be
given, but it may depend on the solution itself as well. To account
for various leak-off regimes, we accept the
following behaviour of $q_l$:
\begin{equation}
\label{q_l} q_l(t,x)=Q_l(t)(l(t)-x)^{\eta},\quad \mbox{for} \quad x
\rightarrow l(t),
\end{equation}
for some constant $\eta\ge-1/2$. Note that the case $\eta=-1/2$
corresponds to the Carter law \cite{Carter}, while $\eta\geq1/3$
guarantees that the leak-off vanishes near the crack tip as
fast as the crack opening at least (see for details \cite{Kusmierczyk}).

The group of fluid equations is to be supplemented by the relation describing deformation of the rock under applied hydraulic pressure. In the case of the PKN model, a
linear relationship between the net fluid pressure and crack
opening is in use:
\begin{equation}
\label{elastic_1} p=kw,
\end{equation}
where a known proportionality coefficient $k=\frac{2}{\pi
h}\frac{E}{1-\nu^2}$ is found from the solution of a plane strain
elasticity problem  \cite{Nordgren} for an elliptical crack of
height $h$. $E$ and $\nu$ are the elasticity modulus and Poisson's
ratio, respectively.

The above equations are equipped with the boundary condition at a crack mouth ($x=0$) determining the injection flux rate:
\begin{equation}
\label{BC_1} -\frac{k}{M}\left[w^3\frac{\partial w}{\partial
x}\right]_{x=0}=q_0(t),
\end{equation}
and two boundary conditions at a crack tip:
\begin{equation}
\label{BC_2} w(t,l(t))=0, \quad q(t,l(t))=0.
\end{equation}
In order to define the crack length, $l(t)$, the global
fluid balance equation is usually utilized (see for example
\cite{Adachi_2007})
\begin{equation}
\begin{array}{l}
\label{fluid balance}
\displaystyle
\int_0^{l(t)}[w(t,x)-w(0,x)]dx-\int_{0}^tq_0(t)dt+
\\[3mm]
\displaystyle
\hspace{30mm}\int_0^{l(t)}\int_{0}^tq_l(t,x)dtdx=0.
\end{array}
\end{equation}
Finally, the initial conditions are assumed in the following way:
\begin{equation}
\label{IC_PKN}
w(0,x)=0, \quad l(0)=0.
\end{equation}

System (\ref{continuos_1}) -- (\ref{IC_PKN}) constitutes the classic formulation of the PKN problem.
It was shown in \cite{Kemp} and \cite{Garagash et al} that the asymptotic behaviours of
 $w$ and  $q$ near the crack tip
are interrelated, and the first term of
the expansion for the crack opening may be written as:
\begin{equation}
\label{As_1} w(t,x)\sim w_0(t)\left(l(t)-x\right)^{\alpha},\quad
\mbox{as}\quad x\to l(t).
\end{equation}
For the classic PKN model the exponent $\alpha=1/3$ was found in
\cite{Kemp}. Thus condition (\ref{BC_2})$_2$ is always
satisfied as it follows from (\ref{Poiseulle_1}) and (\ref{elastic_1}). As a result, the model does not account for the standard
stress singularity of fracture mechanics at the crack tip, and thus
is relevant for the so-called zero toughness regime (see
e.g.\cite{Adachi_2002}).

\vspace{3mm}
 {\sc Remark 1.} Despite that zero crack
opening and length are considered as the initial conditions, all authors begin their studies from the asymptotic model for the
small time. With the assumption of zero leak-off term in the
continuity equation and constant $q_0$, the problem
is reduced to a self-similar formulation. The full numerical
analysis is then continued by taking the similarity solution as the
initial state. This effectively means that the initial conditions
\eqref{IC_PKN} can be replaced by the non-zero crack opening
\begin{equation}
\label{IC}
l(0)=l_{\diamond},\quad
 w(0,x)=w_{\diamond}(x),\quad x\in(0,l_{\diamond}).
\end{equation}

In this paper, the modified formulation of the PKN model is
considered, following the recent advance in the area of numerical
modelling \cite{Linkov_1,Linkov_2,Linkov_3,M_W_L}. Thus, to trace
the fracture front we use the so-called speed equation,
instead of the fluid balance relationship \eqref{fluid balance}:
\begin{equation}
\label{SE} \frac{dl}{dt}=v_0(t)=\frac{q}{w}\big |_{x=l(t)}.
\end{equation}
The speed equation assumes that the fracture tip coincides with the
fluid front, which excludes the presence of a lag or an invasive
zone ahead of the fracture tip. Originally it was introduced by Kemp
\cite{Kemp} and has been recently revisited by Linkov
\cite{Linkov_1,Linkov_2,Linkov_3}.

Note that, on substitution of equations \eqref{Poiseulle_1}, \eqref{elastic_1} and
\eqref{As_1} into \eqref{SE}, one obtains a relationship between the crack
propagation speed and the multiplier of the leading term of the crack
opening asymptotic expansion \eqref{As_1}:
\begin{equation}
\label{SE_1} \frac{dl}{dt}=\frac{k}{3M}w_0^3(t).
\end{equation}
\emph{This implies that the quality of the
numerical estimation of $w_0$ (see estimate (\ref{As_1})) should be vital for the accuracy of
computations.}

By  substituting the Poiseulle equation (\ref{Poiseulle_1}) into the
continuity equation (\ref{continuos_1}) one obtains a lubrication
(Rey\-nolds) equation for the considered problem, where the net fluid
pressure function $p(t,x)$ is eliminated:
\begin{equation}
\label{stiff_1}  \frac{\partial w}{\partial t}
-\frac{k}{M}\frac{\partial}{\partial x}\Big(w^3\frac{\partial
w}{\partial x}\Big)+q_l=0,\quad t \geq t_0,\quad 0 \leq x \leq l(t).
\end{equation}

In this way the modified formulation of the PKN model
includes: i) the Reynolds equation \eqref{stiff_1}; ii)
the boundary conditions \eqref{BC_1} -- \eqref{BC_2}$_1$; iii) the
asymptotics \eqref{As_1}; iv) the initial conditions \eqref{IC}; v)
the speed equation in the form \eqref{SE_1}.

The paper is organized as follows: in the next subsection we present
the normalized formulation of the problem. Then, two types of
self-similar solutions for the PKN model are discussed. These
solutions are used in section 2 to investigate a
numerical algorithm for a time independent
variant of the problem. In section 3, the algorithm is modified to
tackle the transient regime. Two alternative integral
solvers are developed and their performances and applicability
are examined. Section 4 contains the final conclusions.

\subsection{Normalized formulation.}

Following \cite{Kusmierczyk}, we normalize the problem by introducing dimensionless variables:
\[
\tilde x=\frac{x}{l(t)}, \quad \tilde t = \frac{t}{t_n},\quad
t_n=\frac{M}{kl_{\diamond}},\quad \tilde w_\diamond(\tilde x)=
w_\diamond(x),
\]
\begin{equation}
\label{dimensionless}
\tilde w(\tilde t,\tilde
x)=\frac{w(t,x)}{l_{\diamond}},
\quad L(\tilde
t)=\frac{l(t)}{l_{\diamond}},\quad  l^2_{\diamond}\tilde q_0(\tilde t)= t_nq_0(t),
\end{equation}
\[
\tilde
l_{\diamond}\tilde q_l(\tilde t,\tilde x)=t_nq_l(t,x),\quad
l_\diamond^{2/3}\tilde w_0(\tilde t)/L^{1/3}(\tilde t)=w_0(t),
\]
where $\tilde x\in[0,1]$, $L(0)=1$.

\noindent In the new variables equation (\ref{stiff_1}) reads:
\begin{equation}
\label{stiff_2} \frac{\partial \tilde w}{\partial \tilde t}-\tilde x
\frac{L'}{L}\frac{\partial \tilde w}{\partial \tilde x}
-\frac{1}{L^2(t)}\frac{\partial}{\partial \tilde x}\Big(\tilde
w^3\frac{\partial \tilde w}{\partial \tilde x}\Big)+\tilde
q_l=0,
\end{equation}
\[
\tilde t \geq 0,\quad 0 \leq \tilde x \leq 1.
\]
The boundary conditions (\ref{BC_1}) -- (\ref{BC_2})$_1$ may be rewritten
as:
\begin{equation}
\label{BC_11} -\frac{1}{L(\tilde t)}\left[\tilde w^3\frac{\partial
\tilde w}{\partial \tilde x}\right]_{\tilde x=0}=\tilde q_0(\tilde
t), \quad \tilde w(\tilde t,1)=0.
\end{equation}
 The initial conditions (\ref{IC}) are defined as:
 \begin{equation}
\label{IC_11} L(0)=1
,\quad \tilde w(0,\tilde x)=\tilde w_\diamond(\tilde x),\quad
\tilde x\in[0,1].
\end{equation}
The asymptotic expansion for
crack opening (\ref{As_1}) takes the form:
\begin{equation}
\label{As_11} \tilde w(\tilde t,\tilde x)\sim {\tilde w}_0(\tilde
t)(1-\tilde x)^{1/3},\quad \mbox{for} \quad \tilde x \rightarrow 1.
\end{equation}

For the sake of completeness of the normalization, we also present the global fluid balance equation \eqref{fluid balance},
although it will not be used later on:
\begin{equation}
\begin{array}{l}
\label{fluid_balance_1}
\displaystyle
L(\tilde t)\int_0^{1}\tilde w(\tilde t, x)d
x-\int_0^1 \tilde w(0,x)d x-\int_{0}^{\tilde t} \tilde q_0( t)d
t
\\[3mm]
\displaystyle
+\int_{0}^{\tilde t} L(t)\int_0^{1}\tilde q_l( t, x)d x d t=0.
\end{array}
\end{equation}
Finally, the transformation of the speed equation \eqref{SE_1}
yields:
\begin{equation}
\label{SE_2} \frac{d}{dt}L(\tilde t)=V_0(\tilde t)=\frac{1}{3L(\tilde t)}\tilde
w_0^3(\tilde t).
\end{equation}
As shown in \cite{M_W_L}, equation
(\ref{SE_2}) is convenient to trace the fracture front when standard
ODE solvers are in use for the dynamic system describing the
problem. On the other hand, the crack length can be computed from
(\ref{SE_2}) by direct integration to give:
\begin{equation}
\label{CL}
 L(\tilde
t)=\sqrt{1+\frac{2}{3}\int_{0}^{\tilde t}\tilde
w_0^3(\tau)d\tau},
\end{equation}
which is useful when an implicit method (for example Crank-Nicolson
scheme) is utilised (see \cite{M_W_L}).

On substitution of (\ref{SE_2}) into (\ref{stiff_2}) one can rewrite
the later to obtain:
\begin{equation}
\label{Rey} 3L^2\left(\frac{\partial \tilde w}{\partial \tilde
t}+\tilde q_l(\tilde t,\tilde x)\right)=\tilde x \tilde
w_0^3\frac{\partial \tilde w}{\partial \tilde
x}+3\frac{\partial}{\partial \tilde x}\left(\tilde w^3\frac{\partial
\tilde w}{\partial \tilde x}\right).
\end{equation}
In the following, equation
\eqref{Rey} will be used as a basic relation to formulate our
integral solvers.

From now on, for convenience, we shall omit the tilde  symbol in all
quantities. In this way all the notations refer henceforth to the normalized formulation.

\subsection{Self-similar solutions}

Let us assume
\begin{equation}
\label{q_*} q_l(t,x)=\gamma e^{\gamma t}q_l^*(x),
\end{equation}
and look for the similarity solution of the problem in the form:
\begin{equation}
\label{exp}
w(t,x)=u(x)e^{\gamma t},\quad w_0(t)=u_0e^{\gamma t},
\end{equation}
where the asymptotic behaviour (\ref{As_11}) holds true, and $u_0$
is the limiting value of $u$ defined in the same manner as in
estimate (\ref{As_11}). Thus, equation
(\ref{SE_2}) transforms to an identity if one takes:
\begin{equation}
\label{L_exp}
L^2(t)=\frac{2u_0^3}{9\gamma }e^{3\gamma t}.
\end{equation}

On
substitution of \eqref{q_*}, \eqref{exp},  and \eqref{L_exp} into
the equation (\ref{Rey}) one can reduce the latter to the following
ordinary differential equation:
\begin{equation}
\label{Rey1} \beta u_0^3(u+q_l^*)={\cal A}(u),
\end{equation}
with $\beta=2/3$. Here, the nonlinear differential operator ${\cal
A}$ is defined by the right-hand side of
equation (\ref{Rey}) and is equipped with the boundary conditions
\begin{equation}
\label{BC_12} - 3u_0^{-3/2}\left[u^3\frac{d u}{d \tilde
x}\right]_{\tilde x=0}=q_0^*,\quad u(1)=0,
\end{equation}
where we have introduced an auxiliary notation:
\begin{equation}
\label{q_0^*} q_0^*=\sqrt{\frac{2}{\gamma}}e^{-\frac{5\gamma
t}{2}}q_0(t).
\end{equation}
If $q_0^*$ is constant, then  equation
(\ref{Rey1}) together with the boundary conditions (\ref{BC_12}) do
not depend on time and constitute a boundary value problem (BVP) degenerated at point $x=1$.
Indeed, the nonlinear coefficient in front of the second
order term of the differential operator vanishes at the point $x=1$ in view
of the boundary condition (\ref{BC_12})$_2$. This BVP is in fact a self-similar formulation of the original
problem with specific, given leak-off regime and the inlet flux.

Other class of similarity solutions can be found, for some $a\ge0$,
by assuming:
\[
q_l(t,x)=\gamma(t+a)^{\gamma-1} q_l^*(x),\quad w(t, x)=(t+a)^{\gamma}u( x),
\]
\begin{equation}
\label{power}
w_0(t)=u_0(t+a)^{\gamma},
\end{equation}
\begin{equation}
\label{L_power}
L^2(t)=\frac{2u_0^3}{3(3\gamma+1)}(a+t)^{3\gamma+1}.
\end{equation}

As a result, one again obtains the BVP (\ref{Rey1})
-- (\ref{BC_12}) with $\beta=2\gamma/(3\gamma+1)$ and
\begin{equation}
q_0^*=\sqrt{\frac{6}{(3\gamma+1)}}(a+t)^{\frac{1-5\gamma}{2}}q_0(t).
\end{equation}
Thus, for $\gamma=1/5$ the self-similar solution corresponds to the
constant injection flux rate, while the crack propagation speed
decreases with time as $L'(t)=O(t^{-1/5})$ for $t\to\infty$. If, however, one takes $\gamma=1/3$,
the crack propagation speed
is constant and the injection flux rate increases with time: $q_0(t)=O(t^{1/3})$ for $t\to\infty$.

Note that self-similar solutions do not necessarily satisfy the initial
conditions (\ref{IC_11}) as the normalised initial crack lengths are:
\[
L(0)=\sqrt{\frac{2u_0^3}{9\gamma }},\quad
L(0)=\sqrt{\frac{2u_0^3a^{3\gamma+1}}{3(3\gamma+1)}},
\]
for the first and the second type, respectively.

{\sc Remark 2.} As one can see the second type of similarity
solution has a physical sense for any $-1/3<\gamma<\infty$ and thus,
can be used to model three different transient regimes of the crack
evolution: crack acceleration ($\gamma>1/3$), crack deceleration
($\gamma<1/3$), and a steady-state propagation of the fracture
($\gamma=1/3$). The first type of solution possesses a physical
interpretation only for positive values of $\gamma$, which restricts
its application to the cases of accelerating crack.

The self-similar solutions
formulated above are used in the following sections to
analyse computational accuracy provided by the developed solvers.
So far to this end, the asymptotic models have been usually
employed \cite{Nordgren,Kemp,Kovalyshen,Linkov_3}.
However, all of them are restricted to
the case of a constant influx, $q_0$.

\section{Numerical solution of the self-similar problem}

In this section we will formulate an algorithm of the solution for the
self-similar problem defined by equation (\ref{Rey1}) and the
boundary conditions (\ref{BC_12}). The following representation of
the sought function $u(x)$ will be accepted:
\begin{equation}
\label{u_rep} u(x)=u_0(1-x)^{1/3}+\Delta u(x).
\end{equation}
It results from the asymptotic
behaviour (\ref{As_11}) and $\Delta u(x)$ $=O((1-x)^{\zeta})$ for $x
\rightarrow 0$. Parameter $\zeta>1/3$ depends strongly on the behavior
of the leak-off function $q_l$ near the crack tip. In particular,
when $q_l$ vanishes near the crack tip in the same manner as the
solution, or faster, ($\eta\ge1/3$) then $\zeta=4/3$.
One can show that (compare \eqref{q_l})
\begin{equation}
\label{eta}
\zeta=\min\{4/3,1+\eta\}\geq 1/2,
\end{equation}
see also \cite{Kusmierczyk} for details.

\subsection{Integral solver for the self-similar problem.}
\label{subsection_SS_algorithm}

Below we present an algorithm to solve equation (\ref{Rey1}) by
numerical inversion of the operator $\cal A$. Exploiting the
solution representation  (\ref{u_rep}), the inverse operator
$\cal{A}$$ ^{-1}$ defines both components: $u_0$ and $\Delta u$. To
derive $\cal A$$^{-1}$, we integrate the equation (\ref{Rey1})
twice over the interval $[x,1]$ taking into account the boundary
condition (\ref{BC_12})$_2$. Then, after simple transformations one
obtains:
\begin{equation}
\label{dm1}
\begin{array}{l}
3u_0^3 (1-x)\Delta u=-\frac{3}{4} [ 6u_0^2(1-x)^{2/3} (\Delta
u)^2+
\\[1mm]
\displaystyle
4u_0(1-x)^{1/3} (\Delta u)^3+(\Delta u)^4]
+u_0^3\int_{x}^1 \Delta u d\xi+
\\[1mm]
\displaystyle
2u_0^3\int_{x}^1 (\xi-x)u d\xi-(1-x)u_0^3\int_{x}^1 u d\xi+
\\[1mm]
\displaystyle
\beta
u_0^3 \int_{x}^1 (\xi - x)(u+q_l^*)d\xi.
\end{array}
\end{equation}
In short, the latter can be symbolically written in the
compact form:
\begin{equation} \label{invA_1}
\Delta u=G_1(\beta,u_0, \Delta u)+G_2(\beta,u_0, \Delta u,q_l^*),
\end{equation}
where the operators involved in the right-hand side are defined as
follows:
\begin{equation}
\begin{array}{l}
\label{G_1} 3u_0^3(1-x)G_1=-\frac{3}{4}\big[ 6(u_0^2(1-x)^{2/3} (\Delta
u)^2 +
\\[1mm]
\displaystyle
4u_0(1-x)^{1/3} (\Delta u)^3+(\Delta u)^4
\big]+
\\[1mm]
\displaystyle
(2+\beta)u_0^3 \int_{x}^1 (\xi-x)\Delta u d\xi,
\end{array}
\end{equation}
\[
3(1-x)G_2=x \int_{x}^1 \Delta u d\xi
+\frac{3}{28}(3\beta-1 )u_0(1-x)^{7/3}+
\]
\begin{equation}
\label{G_2}
\beta\int_{x}^1 (\xi -
x)q_l^* d\xi.
\end{equation}
One can conclude from (\ref{q_l}) and (\ref{As_11}) that for $x\to1$
\begin{equation}
\label{estimate} G_1=O\left((1-x)^{1+\zeta}\right),\\[1mm]
G_2=O\left((1-x)^{\zeta}\right),
\end{equation}
where $\zeta$ defined in (\ref{eta}).

The relation to compute $u_0$ is
derived by integration of \eqref{Rey1} with respect to $x$ from 0 to
1. Then, taking into account the representation (\ref{u_rep}) and the
boundary condition (\ref{BC_12})$_1$, one can formulate the
condition:
\begin{equation}
\label{G_3} G_3(\beta,u_0,\Delta u, q_0^*)=0,
\end{equation}
where
\[ G_3=\frac{3}{4}(\beta+1)u_0^{5/2}+
\]
\[u_0^{3/2}\left
[(\beta+1)\int_0^1 \Delta u dx+ \beta \int_0^1 q_l^* dx\right
]-q_0^*.
\]
It is easy to prove that for any $q_0^*>0$ and $\beta>-1$ there
exists a unique positive solution $u_0$ of equation (\ref{G_3}),
regardless of the values of the functions $\Delta u$ and $q_l^*$.

The inverse operator ${\cal A} ^{-1}$ is defined, by equations
(\ref{invA_1}) and (\ref{G_3}), as:
\begin{equation}
\label{inverse} [u_0,\Delta u]={\cal A} ^{-1}(\beta,q_l^*,q_0^*).
\end{equation}
Its numerical execution is based on the following iterative
algorithm:
\begin{equation}
\label{Am1}
\begin{cases}  G_3 \left (\beta, u_0^{(i+1)},\Delta
u^{(i)},
q_0^* \right )=0,\\[2mm] \Delta u^{(i+1)}=G_1\left (\beta, u_0^{(i+1)}, \Delta
u^{(i)}\right)+\\[2mm]
\hspace{17mm}G_2\left(\beta, u_0^{(i+1)},\Delta
u^{(i)},q_l^*\right ),
\end{cases}
\end{equation}
where superscripts refer to the consequent iterations. In the first
step we assume that:
\begin{equation}
\label{u_00} \Delta u^{(0)}=0.
\end{equation}
Of course, if any better approximation is available, it can replace
(\ref{u_00}). Note the first relation of (\ref{Am1}) is a
nonlinear algebraic equation that can be solved e.g. by the Newton-Raphson
method, while the second equation is a typical iterative
relationship.

%\subsection{Inverse of the operator $\cal A$. Second algorithm}

\subsection{Numerical examples and discussions}
\label{subsection_numerics_1}

In this section we investigate the performance of the numerical
algorithm (\ref{Am1}) -- (\ref{u_00}). To this end, the power law type self-similar
solution (\ref{power}) is utilized as a benchmark. Apart from the
fact that the equations for other self-similar solution (\ref{exp})
look identically, the value of the parameter $\beta$ appearing in the exponential benchmark is always the same ($\beta=2/3$). Thus the power law type self-similar benchmark gives us
an opportunity to manipulate with the value of $\beta$, which will
be crucial for further implementation of the algorithm to the
transient regime.

Let us utilize the following variant of benchmark solution used previously in \cite{M_W_L}
(compare \eqref{u_bench} in the Appendix A):
\begin{equation}
\label{s_1}
u(x)=(1-x)^{1/3}(1+s(x)),
\end{equation}
\[
s(x)=-\frac{1}{8e}\left(\frac{1}{3}-\frac{2\gamma}{3\gamma+1}\right)(1-x) +0.05(1-x)^2,
\]
which yields $\eta=4/3$.

\noindent Computations are carried out for different
meshes based on $N+1$ nodal points:
\begin{equation}
\label{mesh2} x_j^{(\varrho)}=1-\Big(1-\frac{j}{N} \Big)^\varrho,
\quad j=0,1,...,N.
\end{equation}
When one sets $\varrho=1$, the points of spatial mesh are uniformly
distributed over the whole interval - this mesh will be called the
\emph{uniform mesh}. By taking $\varrho>1$ we obtain increased mesh
density while approaching the end point $x=1$ (the larger
$\varrho$ the greater mesh density near the crack tip) - this mesh
will be refereed to as the \emph{non-uniform mesh}. In our
computations we will use $\varrho=3$. This choice  is motivated by
the following reason. To increase the solution accuracy, we compute
integral operators from (\ref{Am1}) using the classic Simpson
quadrature, which gives an error controlled by the fourth derivative
of the integrand. Accounting for the asymptotic behaviour of the
solution, the transformation (\ref{mesh2}), for $\varrho=3$, improves
the smoothness of the integrand with respect to the new independent
variable near the crack tip (replaces the fraction powers function).
In this way, the error of integration can be minimized. We shall
confirm this below in numerical tests.

\noindent
We use two parameters as the measures of computational accuracy:
 (\emph{i}) -- the maximal point-wise
relative error (denoted as: $\delta u$) of the complete solution
$u(x)$ and
(\emph{ii}) -- a relative error (denoted by $\delta u_0$) of the
coefficient $u_0$ defining the leading asymptotic term in
(\ref{u_rep}).

The efficiency of computations will be  assessed by
the number of iterations needed to compute the solution. The
iterative process is stopped in each case when the $L_2$-norm of the
relative difference between two consecutive approximations
becomes smaller than $\epsilon=10^{-10}$.
Finally, we analyze how the solution errors depend on the number of nodal points $N$.

The computations revealed that convergence of the iterative process
(\ref{Am1}) may only be achieved for some range of $\beta$ values.
For the analyzed benchmark it was: $\beta \in[-1.8,4.8]$. In fact,
this interval is wider than one could expect and fully covers any
physically motivated values of $\beta\in(0,2/3)$  following from the
self-similar formulation.
 Interestingly,
there are also solutions for $\beta<-1$ (compare discussions after
equation (\ref{G_3})). This
information shall be used later on, to construct a solver for
transient regimes.

\begin{figure}[h!]
%M/N=1/300

    \hspace{-2mm}
    \\[2mm]\includegraphics [scale=0.40]{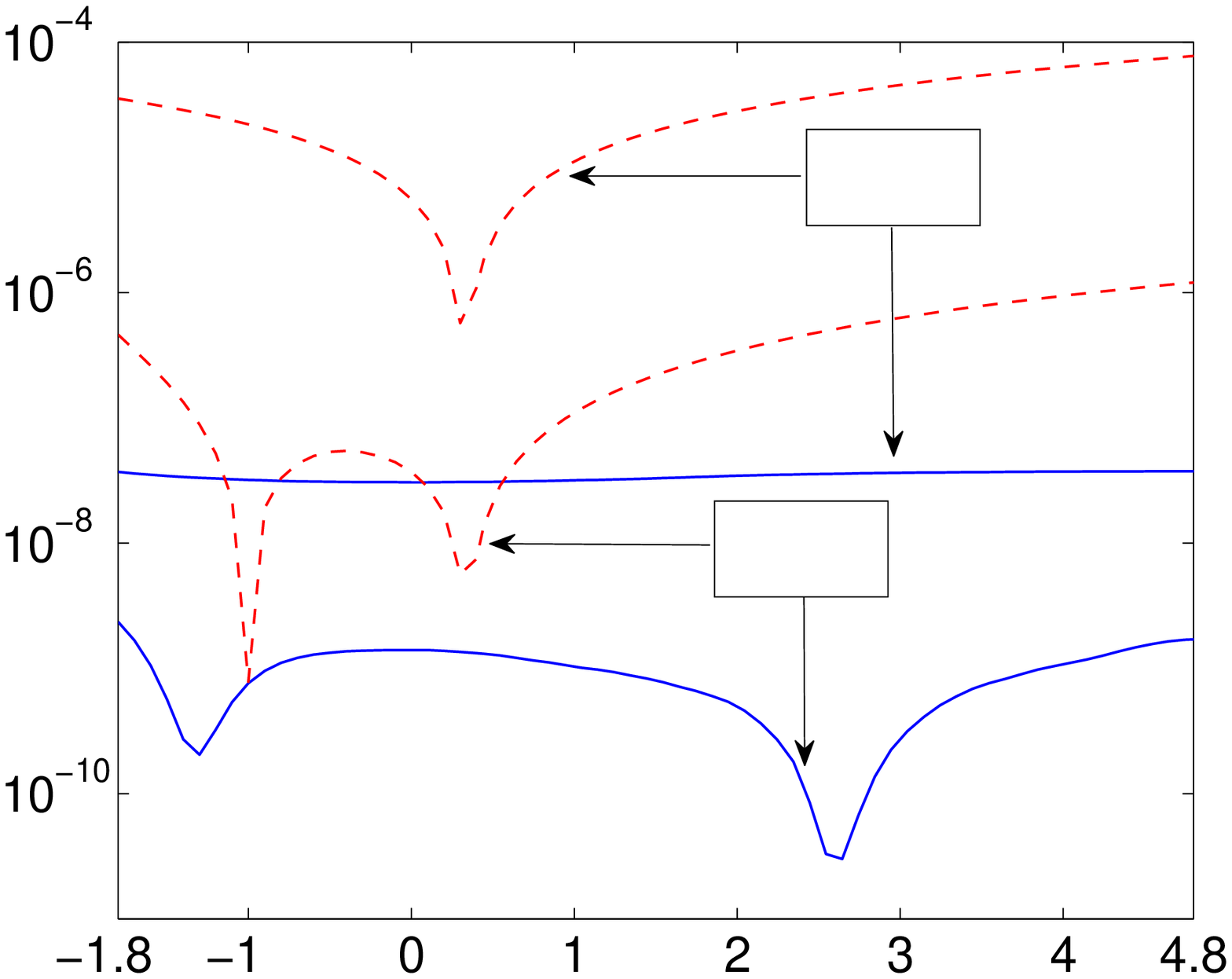}
    \put(-120,-2){$\beta$}
    \put(-90,77){$\delta u_0$}
    \put(-74,137){$\delta u$}
      \put(-225,150){a)}
    \hspace{-2mm}\includegraphics [scale=0.40]{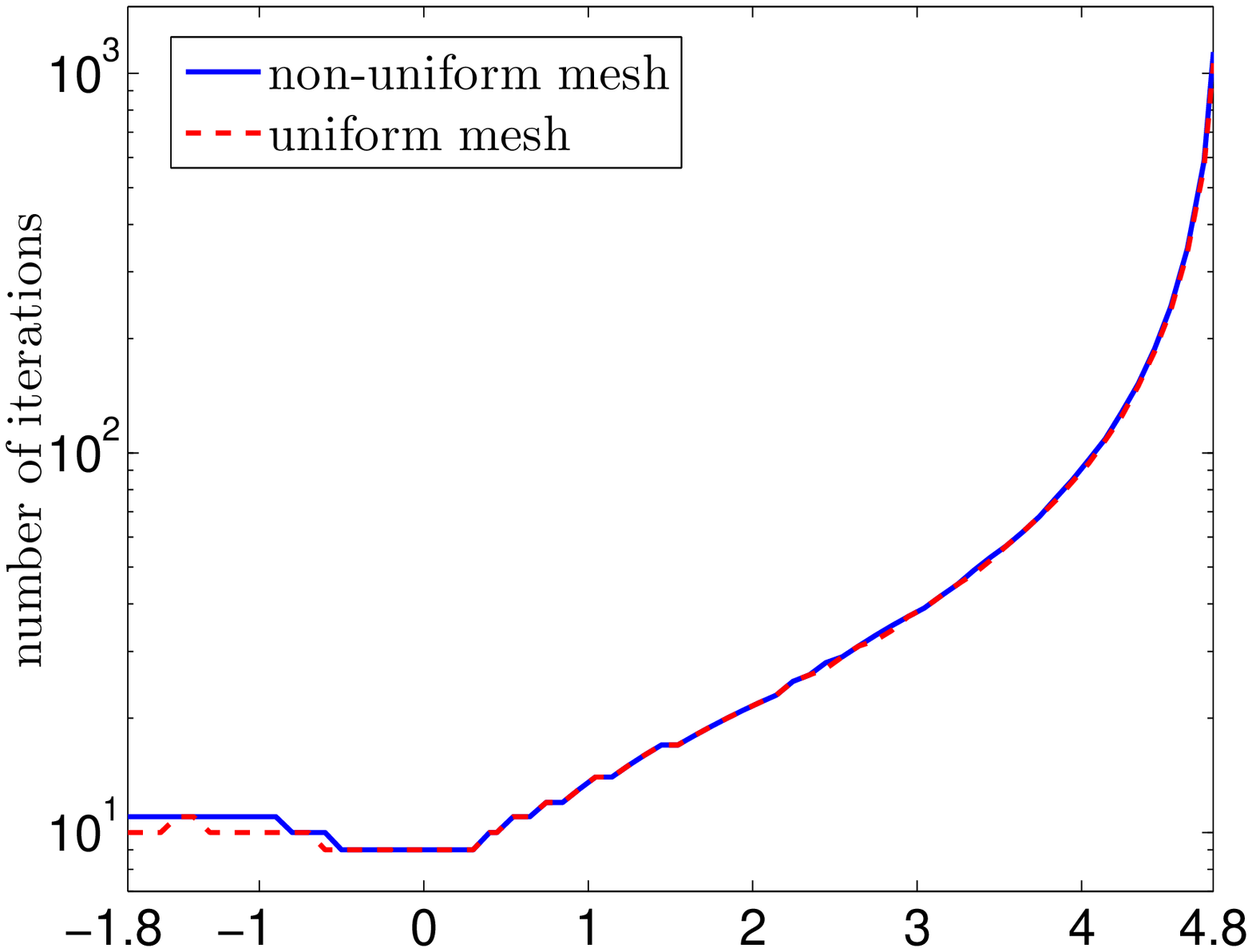}
    \put(-120,-2){$\beta$}
        \put(-225,150){b)}
    \caption{a) The accuracy of computations: maximal relative error of $u(x)$ and $u_0$ as a function of $\beta$;
    b) Number of iterations to reach the final solution as a function of $\beta$. Dashed lines
    refer to the uniform mesh, solid lines to the non-uniform mesh (for $\varrho=3$).  }
\label{del_u}
\end{figure}

\begin{figure}[h!]
%M/N=1/300

    \hspace{-2mm}
    \\[2mm]\includegraphics [scale=0.41]{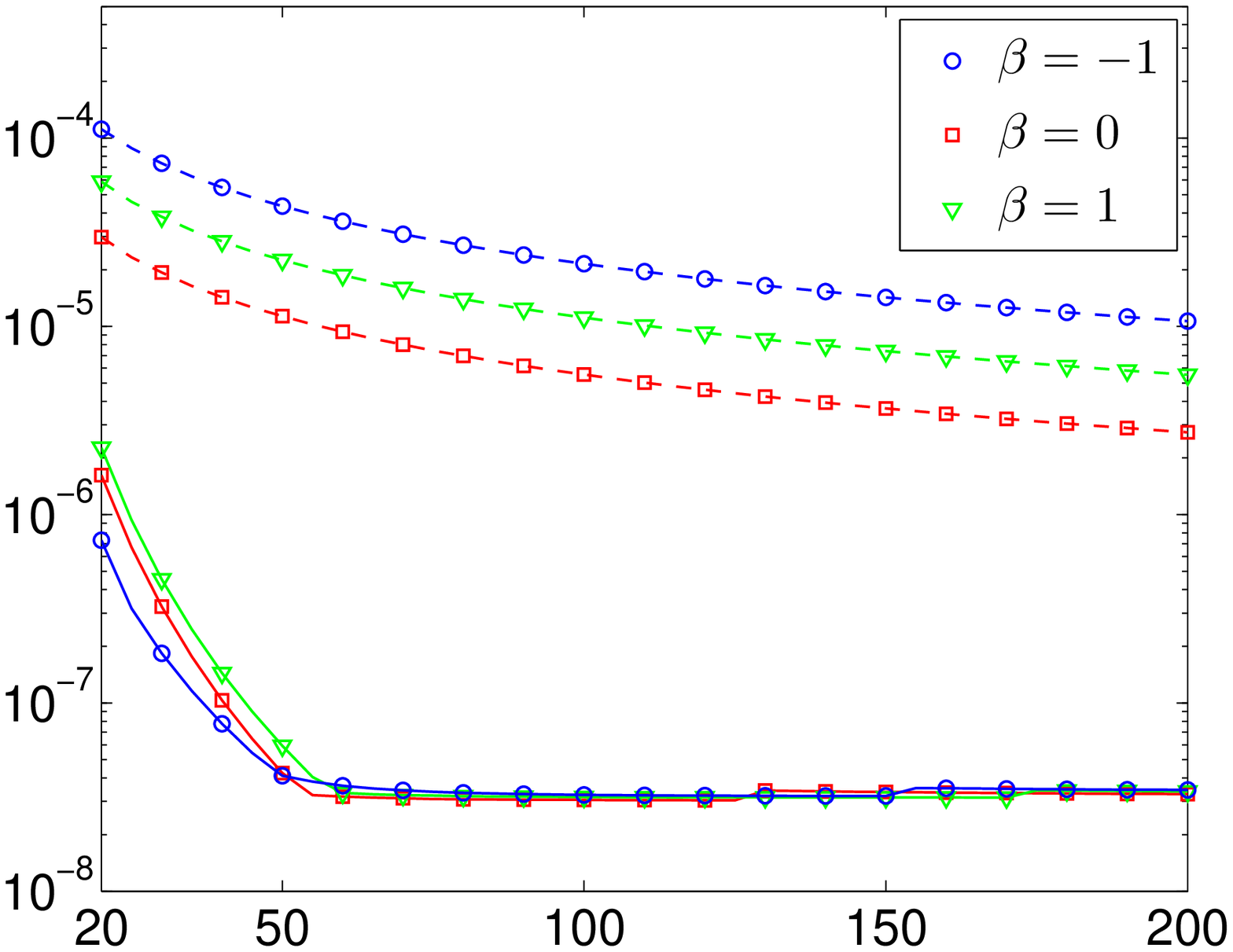}
    \put(-120,-2){$N$}
      \put(-230,150){a)}
    \hspace{-4mm}\includegraphics [scale=0.41]{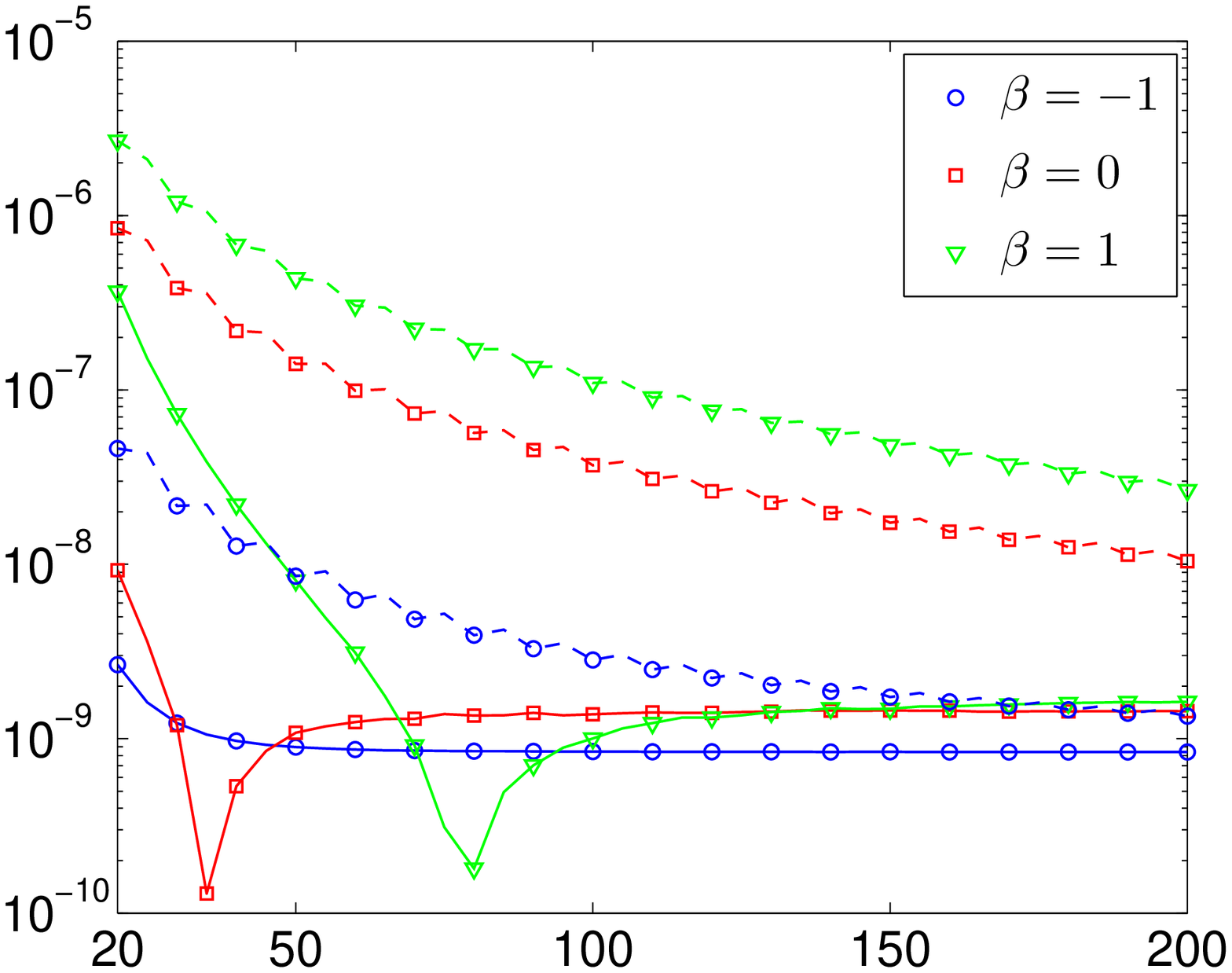}
    \put(-120,-2){$N$}
\put(-230,150){b)}
    \caption{The maximal relative error of solution as a function of the number of nodal points $N$: a) $\delta u$,
    b) $\delta u_0$.  Dashed lines
    refer to the uniform mesh, solid lines to the non-uniform mesh (for $\varrho=3$).}
\label{del_N}
\end{figure}

\begin{figure}[h!]
\center
    %\hspace{-2mm}
    \includegraphics [scale=0.41]{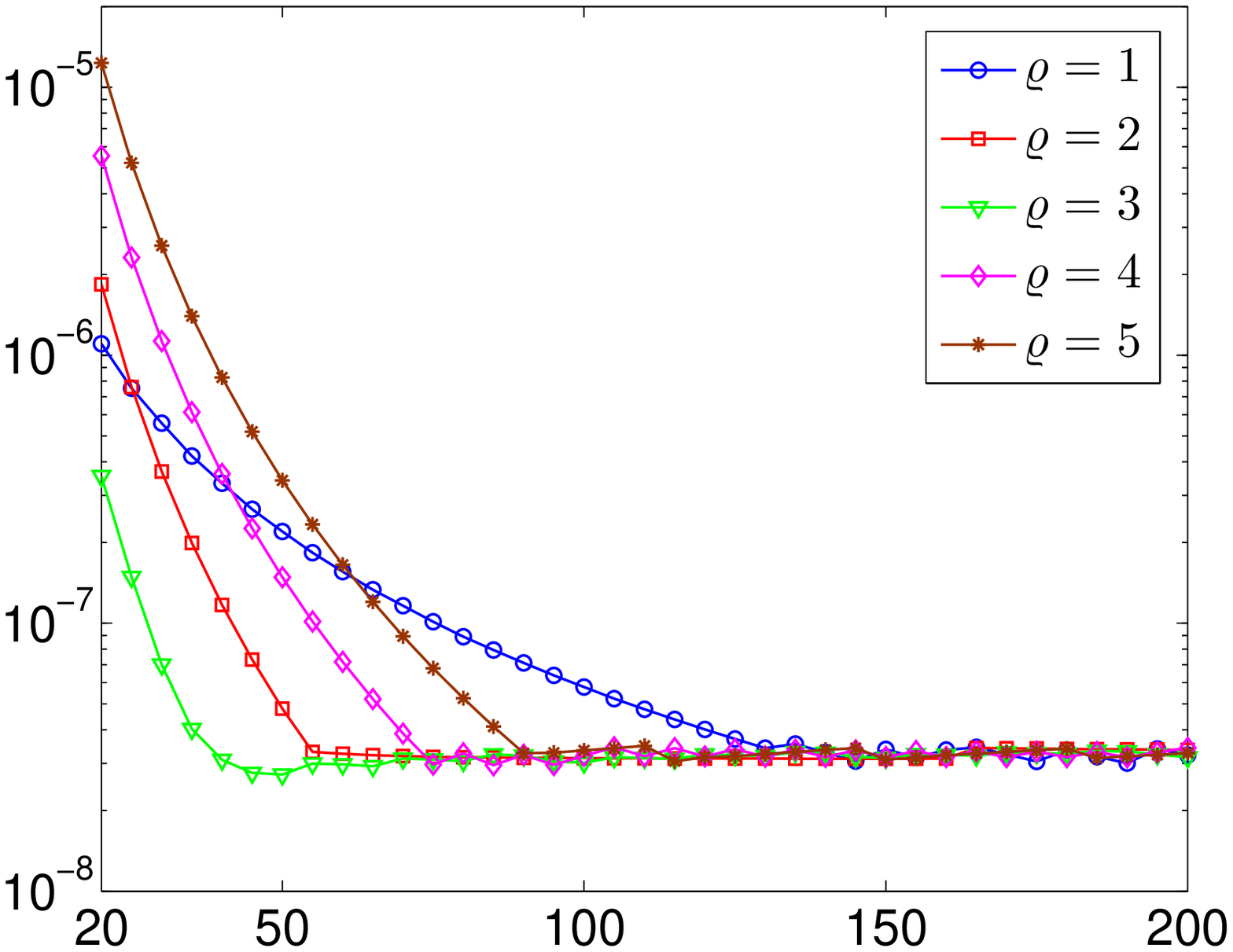}
    \put(-120,-2){$N$}
    \put(-240,90){$\delta u$}

    \caption{The maximal relative error of solution $\delta u_0$ for various spatial meshes. Computations were done for $\beta =1/3$.}
\label{del_N_n}
\end{figure}

In Fig.~\ref{del_u}a) the values of $\delta u$ and $\delta u_0$
obtained for the analyzed benchmark (\ref{s_1}) are shown. The computations were
done for the meshes composed of 100 nodes. As can be seen, the
non-uniform mesh  gives at least one order better accuracy of the
solution, for the same number of nodal points. The error of $u_0$ is
lower than the error of the complete solution $u(x)$, as one could
expect, while its distribution is not as smooth as for $\delta
u(x)$. The minimum of $\delta u(x)$ is located near $\beta=1/3$
which corresponds to the
steady-state similarity solution. Interestingly, in the case of the
uniform mesh the aforementioned minimum is deeper and sharper
than for the non-uniform one.

Fig.~\ref{del_u}b) depicts the number of iterations needed to obtain
the final solution for different values of $\beta$. It shows that
the convergence rate almost does not depend on the type of mesh
chosen. The best efficiency of computations appears for
approximately $|\beta|<1/2$. Note that this interval corresponds to
the values of $\beta$ which provide the best accuracy of
computations. The time of computations (number of iterations)
increases with $|\beta|$ growth, however this trend does not exhibit
a bilateral symmetry (for negative and positive values of $\beta$).
Moreover, the cost of computations increases with the solution error.

Fig.~\ref{del_N} shows how the accuracy of computations depends on the number of nodal points $N$.
Both chosen meshes are analyzed for three different values of $\beta=
-1,0,1$. It can be seen that the non-uniform mesh
 gives at least two orders of magnitude better accuracy than the regular one. For both types of meshes, the values of  $\delta u_0$ are much lower than respective $\delta u$, however the best result for $\delta u_0$ does not necessarily correspond to the best $\delta u$ (e.g. $\beta=-1$ for the regular mesh). The non-uniform mesh provides much lower sensitivity of solution accuracy to the variation of $\beta$ than the uniform one. Also the maximal level of accuracy is obtained much faster for the non-uniform mesh. In the considered case, it is sufficient to take only 60 nodal points to achieve the maximal possible accuracy.

The last test in this subsection identifies the influence of spatial
meshing on the solution accuracy. To this end, we consider the
following  values of $\varrho=1,2,...,5$ from representation
\eqref{mesh2}. The benchmark taken here accepts $\beta=1/3$, as it
 provides the best accuracy and will be important in next
subsections. For each of  the values of $\varrho$, a characteristic
$\delta u(N)$ was computed. The results are presented in
Fig.~\ref{del_N_n}. It shows that, regardless of the mesh under
consideration (or equivalently, the value of the parameter
$\varrho$), the maximal achievable  accuracy is the same. However,
this ultimate level is reached for different numbers of nodal
points, $N$. The fastest convergence takes place for $\varrho=3$,
which confirms our previous predictions on the optimal choice of the
spatial meshing. The slowest convergence to the saturation level
manifests the uniform mesh ($\varrho=1$).

The overall influence of the value of the parameter $\varrho$ on the
accuracy of computations results from the following trend: the
larger the value of $\varrho$ the lower error near the crack tip and the
greater error near the crack inlet. The optimal
balance between the local errors is observed for $\varrho=3$, which
confirms our predictions. In the following, only the non-uniform
mesh for $\varrho=3$ will be used in computations.

\section{Solution to the transient problem}
\label{section_dynamic}

In this section we adopt the idea of the integral solver,
developed for the self-similar formulation, to the
transient regime. The basic assumptions of the approach remain the same,
however the algorithm has to be modified in some essential aspects.
First of all, one has to build the
mechanism of temporal derivative approximation into the numerical procedure together with necessary measures to stabilize the algorithm. We will propose two
methods of doing this, constructing in fact two different solvers.
The second fundamental difference between the self-similar and
time-dependent formulations is that in the latter case, the crack length
$L(t)$ becomes now an element of the solution, which should be looked for
simultaneously with the crack opening $ w(t,x)$ and the first term
of its asymptotics near the crack tip $w_0(t)$.

The basic system of equations for the transient problem is composed
of: the governing equation (\ref{Rey}), the boundary conditions
(\ref{BC_11}), the initial conditions \eqref{IC_11} and the integral
equation defining the crack length (\ref{CL}).

To avoid using multiple subscripts let us adopt the following manner
of notation:
\begin{equation}
\label{not_1} w (t_j,x)=w(x), \quad  w(t_{j+1},x)=W(x).
\end{equation}
Consequently, the asymptotic representations of the solution near the crack tip read:
\begin{equation}
\begin{array}{l}
\label{not_2}  w(x)=w_0(1- x)^{1/3}+\Delta w,
\\[2mm]
\displaystyle
 W(x)=W_0(1-x)^{1/3}+\Delta W,
 \end{array}\quad x \rightarrow 1,
\end{equation}
where $\zeta$ is defined in (\ref{eta}) and
\[
\Delta  w=O((1-x)^\zeta), \quad \Delta W=O((1- x)^\zeta), \quad
\mbox{as}\quad  x \rightarrow 1.
\]

\subsection{Solver using self-similar algorithm (\ref{inverse}) -- (\ref{Am1})}
\label{subsection_IS_1}

Below, we show the way to convert  the initial boundary value
problem defined by equations (\ref{BC_11}), (\ref{CL}), (\ref{Rey})
to the form which may be tackled by the integral solver  in the form
(\ref{Am1}) for the self-similar solution. Obviously, system
(\ref{Am1}) is to be supplemented with an additional equation
defining the crack length.

The main idea of the approach is to use the temporal derivative as one of the dependent variables in the numerical procedure.
To achieve this, we compute the derivative of the solution at each time step $t=t_{j+1}$ in the following iterative process:
\begin{equation}
\begin{array}{l}
\label{temp1}
\displaystyle
\frac {\partial W}{\partial t}^{\!(i+1)}=
G_4\Big(\sigma^{(i+1)},W^{(i+1)},\frac {\partial W}{\partial
t}^{\!(i)}\Big)\equiv
\\[3mm]
\displaystyle
\quad \sigma^{(i+1)}\frac {W^{(i+1)}- w}{\Delta t}+\left
(1-\sigma^{(i+1)}\right )\frac {\partial  W}{\partial t}^{\!(i)},
\end{array}
\end{equation}
where superscripts refer to the number of iteration, $\Delta
t=t_{j+1}-t_{j}$ and the values of  $\sigma^{(i+1)}$ are to be defined later.
The first approximation of the temporal derivative is
\begin{equation}
\label{step_1_1}
\frac {\partial W^{(1)}}{\partial
t}=\frac {\partial  w}{\partial t}.
\end{equation}
 Note
that the derivative at initial time $t=0$
can be immediately found from the governing equation (\ref{Rey}) by
substitution of the initial conditions.

Substituting (\ref{temp1}) into (\ref{Rey}) one obtains the problem
(\ref{Rey1}) and (\ref{BC_12}) with respect to the unknown function
$W(x)$ by exploiting  the following convention:
\[
q_l^*=-w +\frac{\Delta t}{\sigma^{(i+1)}}\left
[(1-\sigma^{(i+1)})\frac{\partial W^{(i)}}{\partial t}+ q_l(t, x)
\right],
\]
\begin{equation}
\label{qs_n}
 \beta W_0^3=\frac{3\sigma^{(i+1)}}{\Delta
t}\left(L^{(i+1)}\right)^2 ,
\end{equation}
and
\begin{equation}
\label{n_beta}
 q_0^*(t)=3q_0(t)W_0^{-3/2}L^{(i)}.
\end{equation}
Finally, the crack length from equation (\ref{CL}) can be iterated as:
\begin{equation}
\label{G_4}
L^{(i+1)}=G_5(W_0)\equiv
\sqrt{\big(L^{(i)}\big)^2+\frac{\Delta t}{3}\left(W_0^3+w_0^3\right)}.
\end{equation}
Here the integral from (\ref{CL}) is computed by the trapezoidal
rule, taking into account the information obtained from the previous
time steps ($t\le t_j$).

In this way the governing partial differential
equation and all the boundary conditions have been transformed to
the formulation used previously in the self-similar problem. As a result,
respective integral operators (\ref{dm1}) -- (\ref{G_2}) and
(\ref{G_3}) remain the same with modified arguments
(\ref{qs_n}) and (\ref{n_beta}), where (\ref{G_4}) should be taken into account.

To choose the value of the parameter $\sigma^{(i+1)}$ at each
iterative step, let us recall that the best performance of the
algorithm described in subsection~\ref{subsection_numerics_1}
has been achieved near $\beta=1/3$. Taking this fact into account in
(\ref{qs_n})$_2$, one can choose
\begin{equation}
\label{sigma} \sigma^{(i+1)}=\frac{ \Delta t \left (W_0^{(i+1)}
\right )^3}{9\left(L^{(i+1)}\right)^2}.
\end{equation}

Thus, the inverse operator for $\cal A$, defining the solution of the
transient problem at the time step $t=t_{j+1}$, can be expressed in the
following manner (compare (\ref{inverse})):
\begin{equation}
\label{inverse1} \Big[W_0, L,\Delta W,\frac{\partial W}{\partial t}\Big]=\hat{\cal A}
^{-1}\big(1/3,q_l^*,q_0^*\big).
\end{equation}

The iterative algorithm of the solver can be described by the system:
\begin{equation}
\label{Am2} \begin{cases}  G_3 \left (1/3, W_0^{(i+1)},\Delta
W^{(i)},
q_0^{*} \right )=0,\\[2mm]
\Delta W^{(i+1)}=G_1\left (1/3, W_0^{(i+1)}, \Delta
W^{(i)}\right)+\\[2mm]
\hspace{18mm}
G_2\left(1/3, W_0^{(i+1)},\Delta
W^{(i)},q_l^{*}\right ),\\[2mm] L^{(i+1)}=G_5\left (W_0^{(i+1)} \right ),\\[2mm]
\displaystyle
\frac {\partial W}{\partial
t}^{\!(i+1)}=G_4\left(\sigma^{(i+1)},W^{(i+1)},\frac {\partial W}{\partial
t}^{\!(i)}\right),
\end{cases}
\end{equation}
Note that parameters $q_0^*=q_0^{*(i)}$, $q_l^*=q_l^{*(i+1)}$ and
$\sigma^{(i+1)}$ are also iterated according to (\ref{qs_n})$_1$,
(\ref{n_beta}) and (\ref{sigma}). To finalize the algorithm, it is
enough to define the initial forms of the crack opening $W^{(1)}$
and the crack length $L^{(1)}$. This is achieved by choosing
\begin{equation}
\label{step_1_2} W^{(1)}=w+\frac{\partial w}{\partial t}\Delta
t,\quad L^{(1)}=L(t_i),
\end{equation}
which finishes the description of the computation scheme for the
fixed time step $\Delta t=t_{j+1}-t_j$. Let us recall that the
temporal derivative of solution at time $t=0$ is computed using the initial conditions (\ref{IC_11}),
while for any next $t=t_j$ we utilize representation (\ref{temp1}).

We would like to underline here the fact that as the output of the
proposed algorithm, one obtains not only the solution of the transient
problem, that is, the crack length $L(t)$ and the crack opening $w(t,x)$, but also the temporal derivative of the latter,
$w_t'(t,x)$ and the crack propagation speed, $V_0(t)$, from (\ref{SE_2}).

\subsection{Solver
based on a modified self-similar algorithm}
\label{subsection_IS_2}

Temporal derivative of the solution can be taken in the following form:
\begin{equation}
\label{temp2} \frac{\partial W}{\partial t}= 2\frac{W- w}{\Delta t}-
\frac{\partial  w}{\partial t},
\end{equation}
which gives the error of approximation of the order $O(\Delta t^2)$. Note that any other
two-points finite difference definition yields only $O(\Delta t)$.

Unfortunately, a direct use of the algorithm formulated in
section \ref{subsection_SS_algorithm} is not, generally speaking,
possible, as it may fail for small time steps. Indeed if one takes
$\sigma=2$ and sufficiently small value of $\Delta t$ in the
representation (\ref{qs_n})$_2$, then the value of the parameter $\beta$ may
be far away from its optimal magnitude $\beta=1/3$. This in turn,
would at least result  in a deterioration of the efficiency of
computations.

On the other hand, it is quite clear that boundary value problem (\ref{Rey1}), (\ref{BC_12}) is solvable for large values of the parameter $\beta$. Indeed,
performing asymptotic analysis, one shows that the solution can be represented in the form:
\begin{equation}
\label{sol_as} u(x)=-q_l^*(x)+b_0(x)+b_1(x),\quad x\in(0,1),
\end{equation}
where $b_0$ and $b_1$ refer to two boundary layers, accounting for the boundary conditions (\ref{BC_12})$_1$ and (\ref{BC_12})$_2$, respectively.

Thus the algorithm should be modified to be able to deal with large
values of $\beta$. To achieve this goal, let us look at the original
representation (\ref{dm1}), where only the  last term in the
right-hand side depends on $\beta$. This term violates the
convergence of the iterative process for large $\beta$. To prevent
this from happening, at each iteration we supplement the term in
question with an auxiliary so-called 'viscous' term in the form
$\mathcal{V}(x)= \beta u_0^3(1-x)(C_0+C_1(1-x))$, where the
constants, $C_0,C_1$, are computed by comparing both the original and the viscous
terms. In this way one can construct a \emph{modified} algorithm,
 schematically represented in the following manner
\begin{equation}
\label{Am3} \begin{cases}  G_3 \left (\beta^{(i)}, W_0^{(i+1)},\Delta
W^{(i)},
q_0^{*} \right )=0,\\[2mm]
\mathcal{V}^{(i+1)}=\mathcal{V}(W^{(i)},\hat q_l^{*(i+1)}),
\\[2mm]
\Delta W^{(i+1)}=\hat G_{12}\left(\beta^{(i+1)}, W_0^{(i+1)},\Delta
W^{(i)},\hat q_l^{*},\mathcal{V}^{(i+1)}\right),\\[2mm]
L^{(i+1)}=G_5\left (W_0^{(i+1)} \right ),\\[2mm]
\displaystyle
\frac{\partial W}{\partial
t}^{\!(i+1)}=2\frac{W^{(i+1)}-w}{\Delta t}-\frac{\partial w}{\partial
t},
\end{cases}
\end{equation}
where
\begin{equation}
\begin{array}{l}
\label{qs_nn} \displaystyle
\hat q_l^*=-w -\frac{\Delta t}{2}\left [\frac{\partial
w}{\partial t}+q_l(t,x) \right],
\\[4mm]
\displaystyle
\beta^{(i)} =\frac{6}{\Delta t
}\left(L^{(i)}\right)^2\left(W_0^{(i)}\right)^{-3}.
\end{array}
\end{equation}
The function $q_0^*(t)$ is defined in the same way as previously (see (\ref{n_beta})). At each time step the initial value of $W$ is taken in the form (\ref{step_1_2}).
We do not show here the precise definition of the  operator $\hat
G_{12}$, as it may be easily derived by merging operators $G_1$ and
$G_2$ with the viscous term $\mathcal{V}(x)$, where respective
constants are computed, for example by the least squares method.

\subsection{Analysis of the algorithms performance.}
\label{subsection_numerics_2}

The aim of this subsection is to analyze and compare the performances of the algorithms formulated in
subsections \ref{subsection_IS_1} and \ref{subsection_IS_2}.
To this end, the benchmark solution used previously for the
self-similar formulation is utilized, for the time dependent term
$ \psi(t)=(1+t)^{\gamma}$ (compare \eqref{u_bench}). First tests are performed for $\gamma=1/5$.

In the following, the notations \emph{solver 1} and \emph{solver 2} are attributed to the
algorithms \eqref{Am2} and \eqref{Am3}, respectively.

Let us analyze the influence
of the spatial mesh density on the accuracy of computations done by
both solvers. For this reason, we start with a single time step (in
our case $\Delta t=10^{-2}$) and carry out the computations for
different numbers of nodal points $N$,
ranging from 10 to 100. The
following parameters are used for the comparison: the maximal
relative error of the crack opening, $\delta w$, the relative error
of the crack length, $\delta L$, and the maximal relative error of
the temporal derivative of the crack opening, $\delta w_t$.

\begin{figure}[h!]
%M/N=1/300

    \hspace{-2mm}
    \\[2mm]
    \includegraphics [scale=0.40]{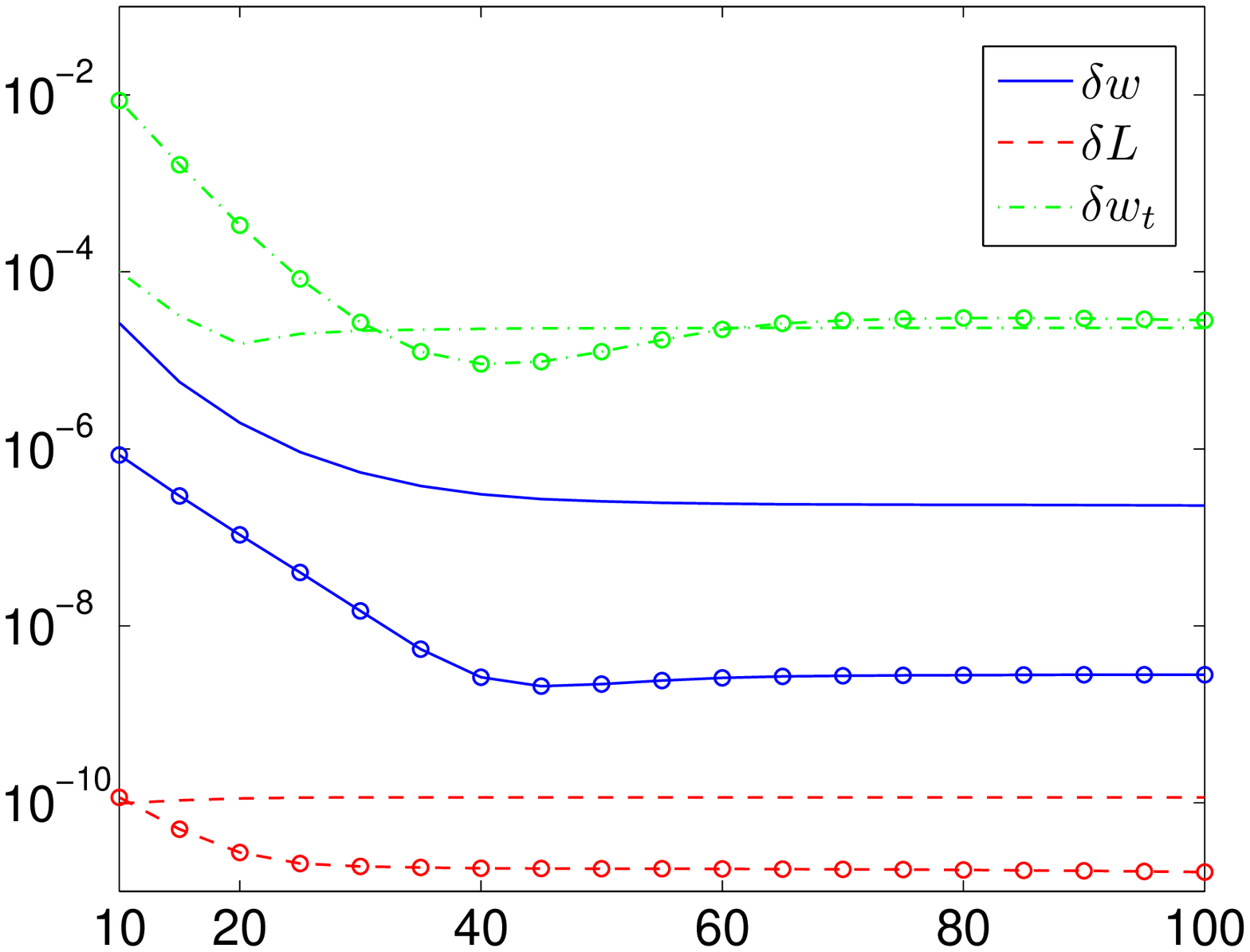}
    \put(-120,-2){$N$}
    \put(-225,150){a)}
    \hspace{2mm}
    \includegraphics [scale=0.40]{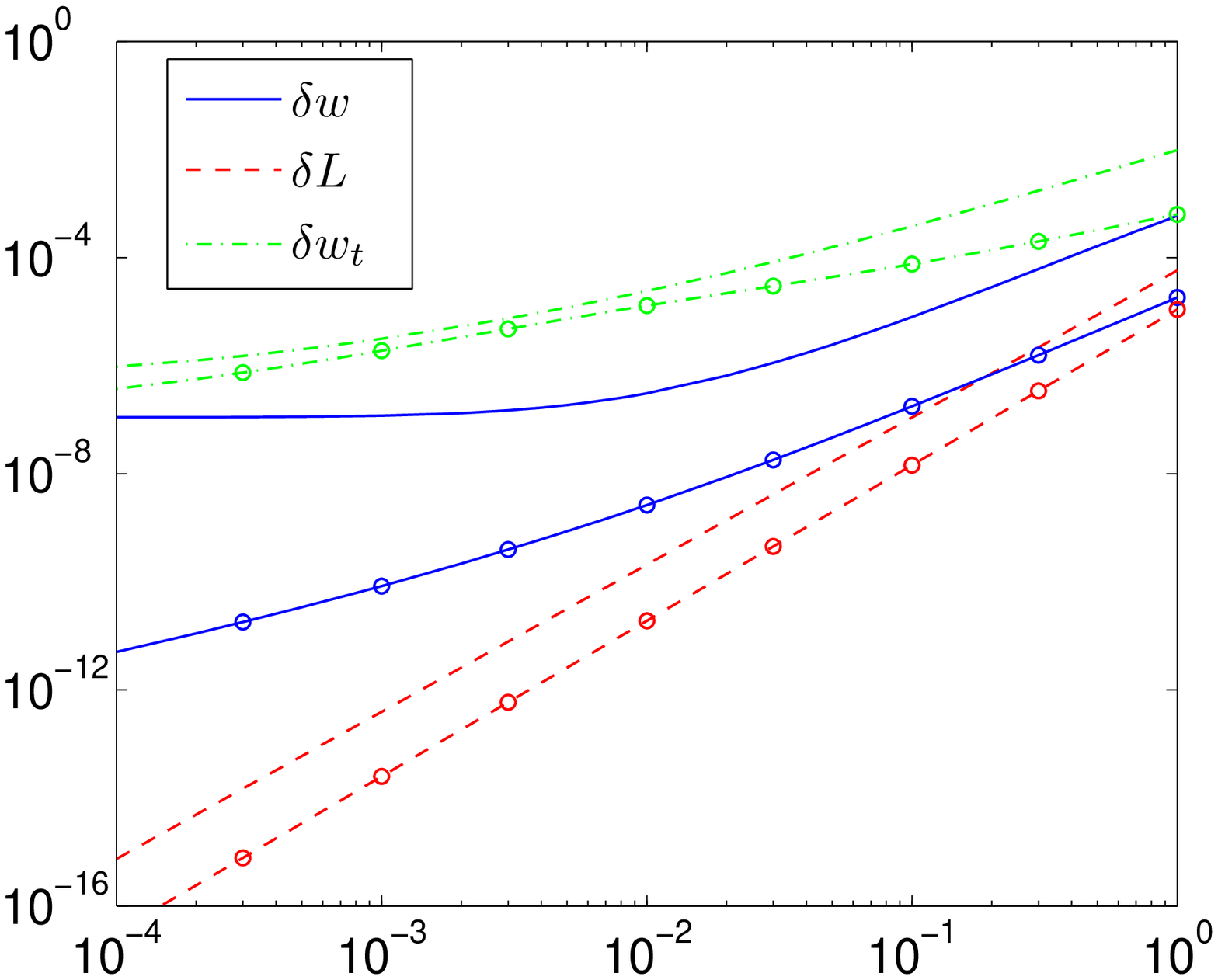}
    \put(-110,-2){$\Delta t$}
    \put(-225,150){b)}

    \caption{The errors of computations: $\delta L, \delta w, \delta w_t$ as  functions of a) spatial mesh density (number of nodal points, $N$) for the fixed time step $\Delta t=10^{-2}$;
    b) time step $\Delta t$ for fixed number of the nodal points $N=40$. Line with markers refer to \emph{solver 2}}
\label{blad_N}
\end{figure}

The results of computations are depicted in Fig.~\ref{blad_N}a). It
shows that regardless of the considered parameter, \emph{solver 2}
always provides  better accuracy (two orders of magnitude for the
crack opening, $w$, and an order for the crack length, $L$). The
only exception is for $\delta w_t$, whose values are of the same
order and \emph{solver 1} may even give a bit lower errors. However,
this does not result in a better accuracy of the two remaining
components of the solution: $w$ and $L$. Note that for the first
time step, the accuracy of computation of the crack length, $L(t)$,
is of two orders of magnitude better than that for the crack
opening, $\delta w$, regardless of the solver type and the number of
the nodal points, $N$.

Similarly to the trend observed in the self-similar formulation,
computational errors stabilize for some critical value $N=N_*(\Delta
t)$. Surprisingly, for the transient regime where the temporal
derivative plays a crucial role, the critical $N_*$ appears to be
even slightly lower than that for similarity solution. Thus, for
both solvers it is sufficient to take only 40 points to achieve the
maximal level of accuracy.

\vspace{2mm}

{\sc Remark 2}. From first glance it may be surprising that the
accuracy of the temporal derivative, $\delta w_t$, is up to four
orders of magnitude worse than that for $\delta w$ in the case of
\emph{solver 2}. However, this fact can be easily explained when one
analyses the  estimation:
 $\delta w_t\approx 2 w\delta w /(\Delta t
w_t)$ which follows from (\ref{temp2}) for small $\Delta t$.
Computing the multiplier in the right hand side of the estimation,
one obtains the respective value of the order of $10^{4}$. Note
that the recalled formula does not account for the error of the
method of $w_t$ approximation itself. In some cases (very small time
steps) this value may be comparable to the former, and thus
essentially influence the overall error.

\vspace{2mm}

{\sc Remark 3}. At the first time step, there also exists a direct
relationship between the error of the crack length, $\delta L$, and
the multiplier of the leading term of the asymptotics (\ref{not_2}),
$\delta W_0$. For small values of $\Delta t$ it reads: $2\delta
L=W_0^3\delta W_0\Delta t$. Having this relation we do not show a
separate analysis for $W_0$ and consequently, $V_0$.

All the results presented above were obtained for a single time step
$\Delta t=10^{-2}$. In order to illustrate the influence of  $\Delta
t$ on the solution accuracy, a number of computations were done in
the interval $\Delta t \in[10^{-4},1]$ for a fixed number of
nodal points, $N=40$. The results describing the solution errors are
shown in  Fig.~\ref{blad_N} b). One can conclude from them that the
errors decrease when reducing the time step (at least in the
analyzed range of parameters). Simultaneously, one can expect that
this tendency should have its own limitation for a fixed number of
the nodal points, $N$, and starting from some small value $\Delta
t(N)$, it will reverse to the opposite. For  both solvers, a fast
decrease of  $\delta L$ is observed as  $\Delta t$  gets smaller
($\delta L\approx a10^{-4}\Delta t^3$ and $a\sim 1$).

Now let us analyze the accuracy of the solution on a given time
interval $t\in[0,t_K]$. As, it follows from the results presented in
\cite{Kusmierczyk} and \cite{M_W_L}, after the initial growth to the
maximal value, the error of computations stabilizes at some level or
even decreases, and thus further extending of the time interval does
not contribute to the deterioration of accuracy. Bearing this in
mind we take here $t_K=100$. First we test how the number of time
steps $K$ affects the accuracy of computations within the same time
stepping strategy. For both solvers we took a fixed number of
spatial mesh points, $N=40$, changing the number of time steps $K$
from 10 to 300. The utilized time stepping strategy was the same as
that  used in \cite{M_W_L} for the Crank-Nicolson scheme.  It
is given by the following equation ($i=1,2,...,K$):
\begin{equation}
\label{t_step}
t_i=(i-1)\delta t +\frac{t_{K}-(K-1)\delta t}{(K-1)^3}(i-1)^3,
\end{equation}
where $\delta t$ is a parameter controlling the first time step.
Note that by increasing the value of $K$ one distributes the time
points near $t=0$ almost uniformly.

\begin{figure}[h!]
\center
    %\hspace{-2mm}
    \includegraphics [scale=0.40]{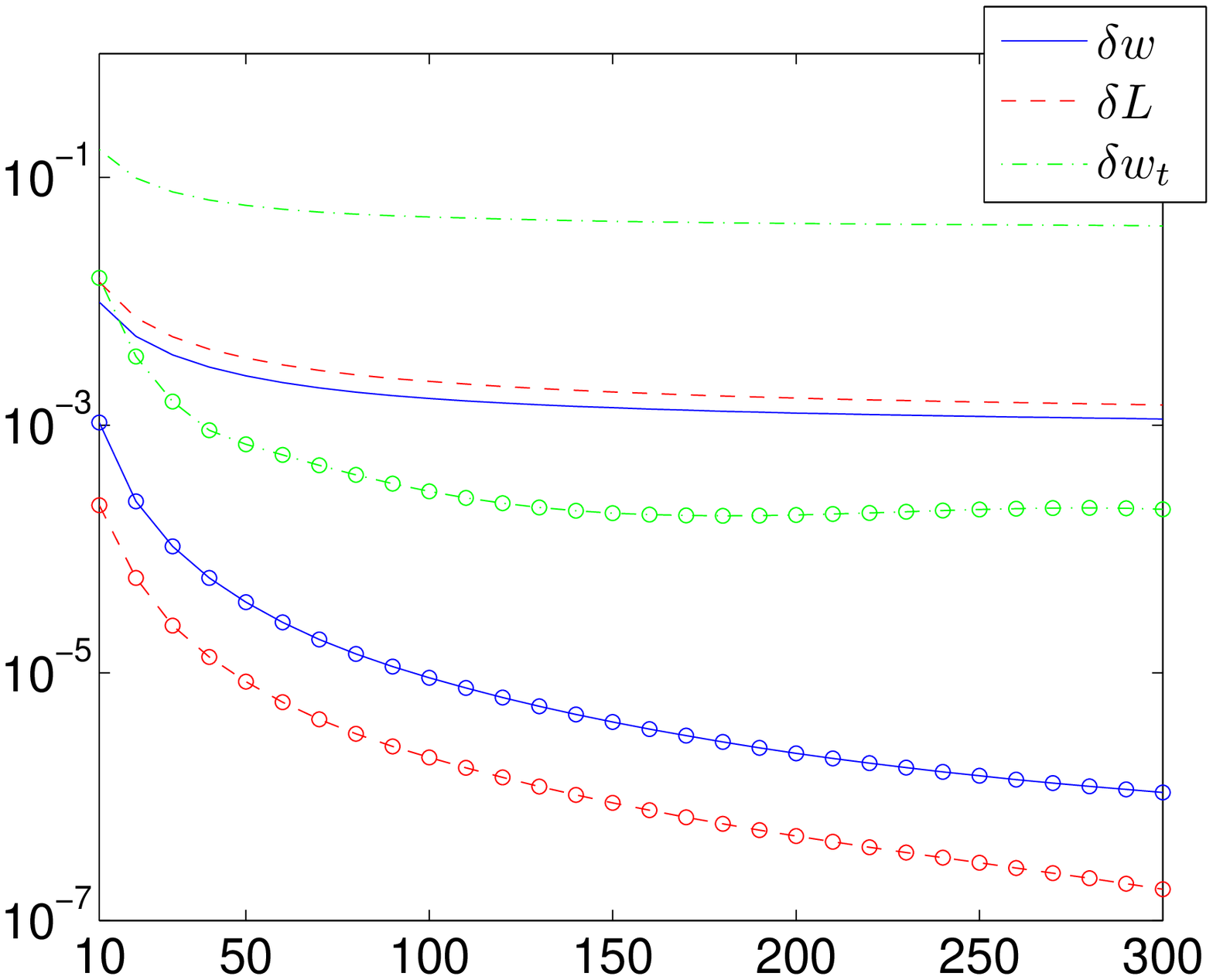}
    \put(-115,-2){$K$}
    \caption{Influence of the number of time steps, $K$, within time stepping strategy (\ref{t_step}), on accuracy of the computations for the fixed non-uniform spatial mesh ($N=40, \varrho=3$).
    Lines without markers correspond to \emph{solver 1}, lines with markers refer to \emph{solver 2}.}
\label{bl_dt}
\end{figure}

The achieved accuracy, as a function of the number of time steps
$K$, is presented in Fig.~\ref{bl_dt}. It shows that for the same
value of $K$ \emph{solver 2} provides much more accurate
results again and this advantage increases with growing $K$. One can
expect, as it follows from the discussions after Fig.~\ref{blad_N}
which refers to a single time step, that each solver gives the
solution errors satisfying the estimation: $\delta L<<\delta
w<<\delta w_t$ for any set of the input parameters. Indeed,
\emph{solver 2} supports this statement as can be seen
Fig.~\ref{bl_dt}, while for \emph{solver 1} surprisingly another
trend is observed: $\delta w<\delta L <<\delta w_t$.

\begin{table*}[t!]
\centering
\begin{tabular}{|c|c|c|c|c|c|c|c|c|}
  \hline
  solver& $N$ & $K$ & $\delta L$ & $\delta w$ &$\delta V_0$&  $\delta w_t$ &$\delta w_t^{({F\!D}_2)}$&$\delta w_t^{({F\!D}_3)}$\\
  \hline  \hline % or \cline{col1-col2} \cline{col3-col4} ...
  {from \cite{M_W_L}}
  & $100$ & $242$ & $6.8\cdot10^{-3}$ & $5.0\cdot 10^{-3}$ &--& -- &$\begin{array}{c}4.1\cdot 10^{-1}\\(1.0\cdot10^{-2})\end{array}$&$\begin{array}{c}4.1\cdot 10^{-1}\\(4.8 \cdot 10^{-3})\end{array}$\\
  \hline
  \emph{solver 1} & $40$ & $30$ & $5.2 \cdot 10^{-3}$ & $3.7\cdot 10^{-3}$ & $7.1\cdot10^{-3}$ &$7.5\cdot 10^{-2}$&$4.6\cdot 10^{-2}$&$2.0\cdot 10^{-2}$\\
  \hline
  \emph{solver 2} & $40$ & $30$ & $2.4 \cdot 10^{-5}$  & $1.1\cdot 10^{-4}$ & $3.2\cdot10^{-4}$  &$1.5\cdot 10^{-3}$&$4.6\cdot 10^{-2}$&$3.4\cdot 10^{-3}$\\
  \hline
  \emph{solver 1} & $5$ & $30$ & $5.2 \cdot 10^{-3}$ & $3.7\cdot 10^{-3}$ & $7.0\cdot10^{-3}$ &$7.5\cdot 10^{-2}$&$4.1\cdot 10^{-2}$&$2.0\cdot 10^{-2}$\\
  \hline
  \emph{solver 2} & $5$ & $30$ & $8.0 \cdot 10^{-5}$  & $5.7\cdot 10^{-4}$ & $6.3\cdot10^{-4}$  &$5.6\cdot 10^{-2}$&$6.3\cdot 10^{-2}$&$6.2\cdot 10^{-2}$\\
  \hline
\end{tabular}
\caption{Comparison of the results obtained for the solver developed in \cite{M_W_L} and two  integral solvers:
\emph{solver 1} and \emph{solver 2}. The same benchmark solution
with the leak-off vanishing near the crack tip is considered.}
\label{table_1}
\end{table*}

\begin{figure}[h!]
\center
    %\hspace{-2mm}
    \includegraphics [scale=0.40]{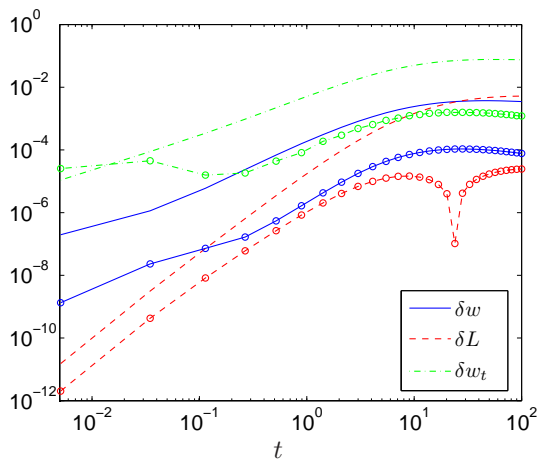}
    \put(-115,-2){$t$}
    \\[1mm]
    \caption{Distribution of the solution errors in time. Line with markers refer to \emph{solver 2}. The time stepping strategy (\ref{t_step}) is used for $K=30$.
    The number of nodal points, $N=40$.}
\label{bl_czas}
\end{figure}

To explain this apparent paradox, let us recall that all errors in
Fig.~\ref{bl_dt} are taken as the maximal values over the time-space
domain.  However, since for the assumed time stepping strategy the
time steps increase with growing time, one can expect (compare the
results in the Fig.~\ref{blad_N}b)) that at some point $\delta L$
may become greater than $\delta w$. Moreover, the error
accumulation in successive time steps may additionally influence the
relationship between $\delta L$ and $\delta w$, especially if they
are close to each other as it is in the case of \emph{solver 1}.
Fig.~\ref{bl_czas} depicts the evolution of computational errors in
time. For both solvers $\delta w$ and $\delta w_t$ reach their
maximal values and stabilize or decrease for $t<t_K$. In the case of
\emph{solver 1} one can observe an intensive error accumulation for
the crack length, $\delta L$. Indeed, it can be seen from the
Fig.~\ref{bl_czas}, that there exists a moment when the error curves
for $\delta L$ and $\delta w$ intersect for the \emph{solver 1}. In
the case of \emph{solver 2}, respective errors of solution ($\delta L$,
$\delta w$ and $\delta w_t$) have essentially different values. As a
result, this effect does not take place in the time interval under
consideration. However, it may be encountered for $t_K>100$, as
suggests the trend observed for $t$ close to $t_K$.

It is interesting to compare the solution errors generated by the two
integral solvers with those obtained  in \cite{M_W_L} for a dynamic
system (DS) approach. In the latter case a standard MATLAB solver,
\emph{ode15s}, was employed. The best results were obtained for a
uniform spatial mesh, the presented data corresponded to $N=100$.
The time stepping strategy chosen automatically by the solver (whose
character is approximated by \eqref{t_step}) accepted 242 time
steps. The maximal relative error of the crack opening, $\delta w$,
and the crack length, $\delta L$, were $5.0\cdot 10^{-3}$ and
$6.8\cdot10^{-3}$ (see Table~\ref{table_1}), respectively.

When analyzing the data collected in Table~\ref{table_1}, one
concludes that it is sufficient for any of the integral solvers to
take only $N=40$ nodal points and $K=30$ time steps, as suggested the previous analysis, to have better
(but comparable -- \emph{solver 1}), or much better (\emph{solver
2}) results. Indeed, the corresponding maximal errors for the
integral solvers are: $\delta w=3.7\cdot 10^{-3}$, $\delta L=5.2
\cdot 10^{-3}$ for \emph{solver 1}, and $\delta w=1.1\cdot
10^{-4}$, $\delta L=2.4 \cdot 10^{-5}$ for \emph{solver 2}. In
other words, the first solver provides the same accuracy for the crack opening and the crack length as the DS
solver using much greater numbers of nodal points and time steps, while
the second one, under the same conditions, improves the results at
least one order of magnitude.
It also shows that the \emph{solver 2} yields one order of magnitude better accuracy of the crack propagation speed,
$V_0$, than the \emph{solver 1}.

The new algorithms allow us to automatically compute the temporal
derivatives in the solution process. The respective errors, $\delta
w_t$, are: $7.5 \cdot 10^{-2}$-- \emph{solver 1} and $1.5 \cdot
10^{-3}$ -- \emph{solver 2}. We decided to compare these figures,
with the ones obtained in postprocessing (here, also the DS
approach was examined). To this end two FD schemes (2-points and
3-points) were used. This time, the corresponding errors, $\delta
w_t$, were: $4.6\cdot 10^{-2}$ and $2.0\cdot 10^{-2}$ for
\emph{solver 1}, $4.6\cdot 10^{-2}$ and $3.4\cdot 10^{-3}$ for
\emph{solver 2} and the same value $4.1\cdot 10^{-1}$ for both
schemes in the case of DS. It is worth mentioning that the values
obtained for DS appeared at the first time step. Then, the errors
decreased with time and  stabilized to give their minimal levels
of $10^{-2}$ and $4.8 \cdot 10^{-3}$, correspondingly.

 As can be seen, the integral solvers give at
least one order of magnitude better accuracy of $w_t$ than the DS. Moreover,
while the postprocessing gives smaller (but comparable) error for
the \emph{solver 1}, \emph{solver 2} returns more accurate values of
$w_t$ than those obtained in the postprocessing, even for the
3-points FD. Finally, apart from the fact that $\delta w$ and
$\delta L$ for the DS and \emph{solver 1} look comparable in values,
the quality of the computation is better for the new solver as is
clear from the postprocessing analysis.

Just for comparison we also present in Table~\ref{table_1} the results obtained for the spatial mesh composed of only five nodal points, $N=5$. It turned out that even for such a drastic reduction  of the mesh density, the solution accuracy for most of the parameters is of the same order as for $N=40$. Interestingly, \emph{solver 1} exhibits almost no sensitivity to this mesh reduction. In fact the distinguishable differences can  be observed only for the \emph{solver 2}.

\begin{figure}[h!]
%M/N=1/300
    \includegraphics [scale=0.40]{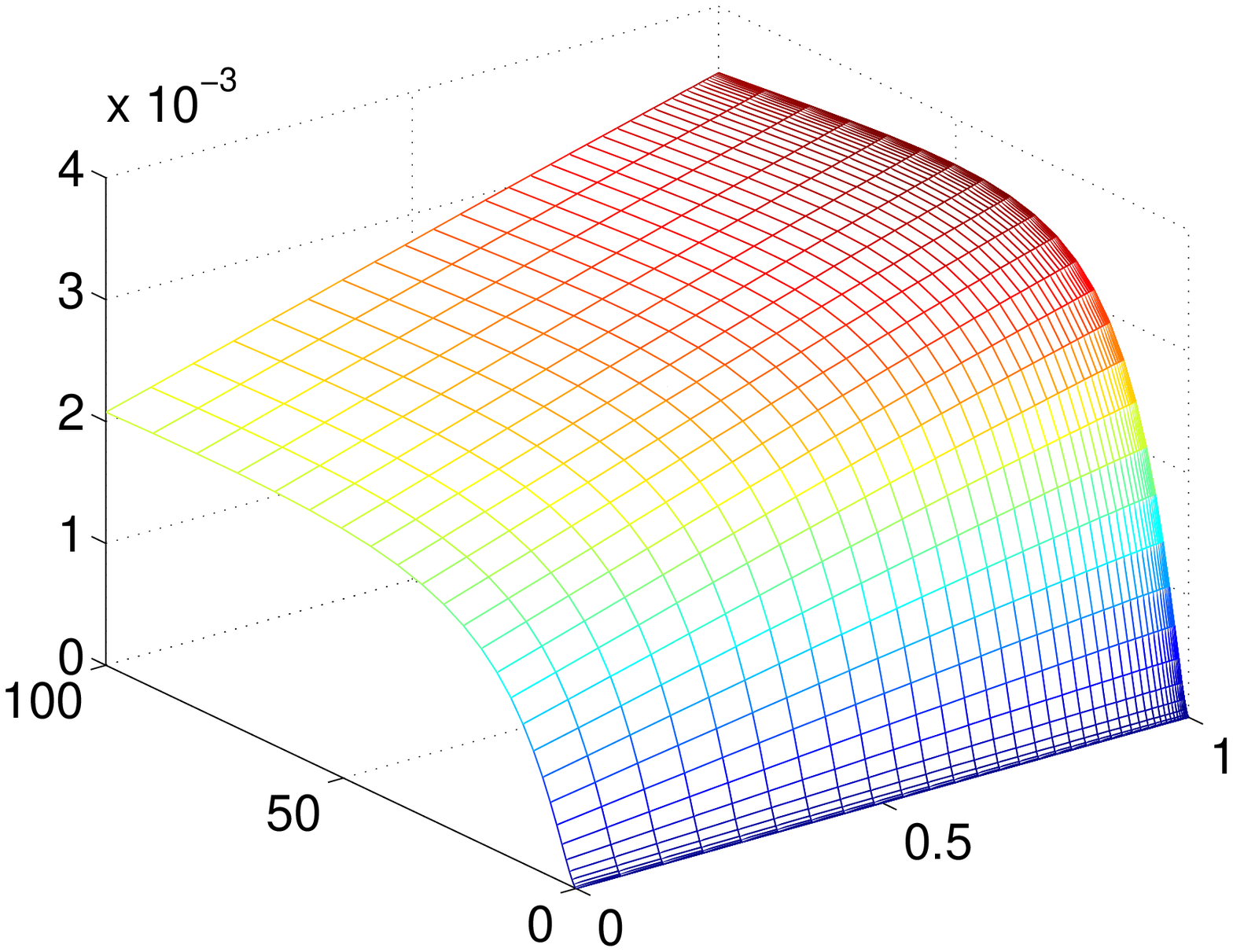}
    \put(-60,10){$x$}
    \put(-180,15){$t$}
    \put(-225,95){$\delta w$}
    \put(-225,150){a)}
    \hspace{-2mm}\includegraphics [scale=0.40]{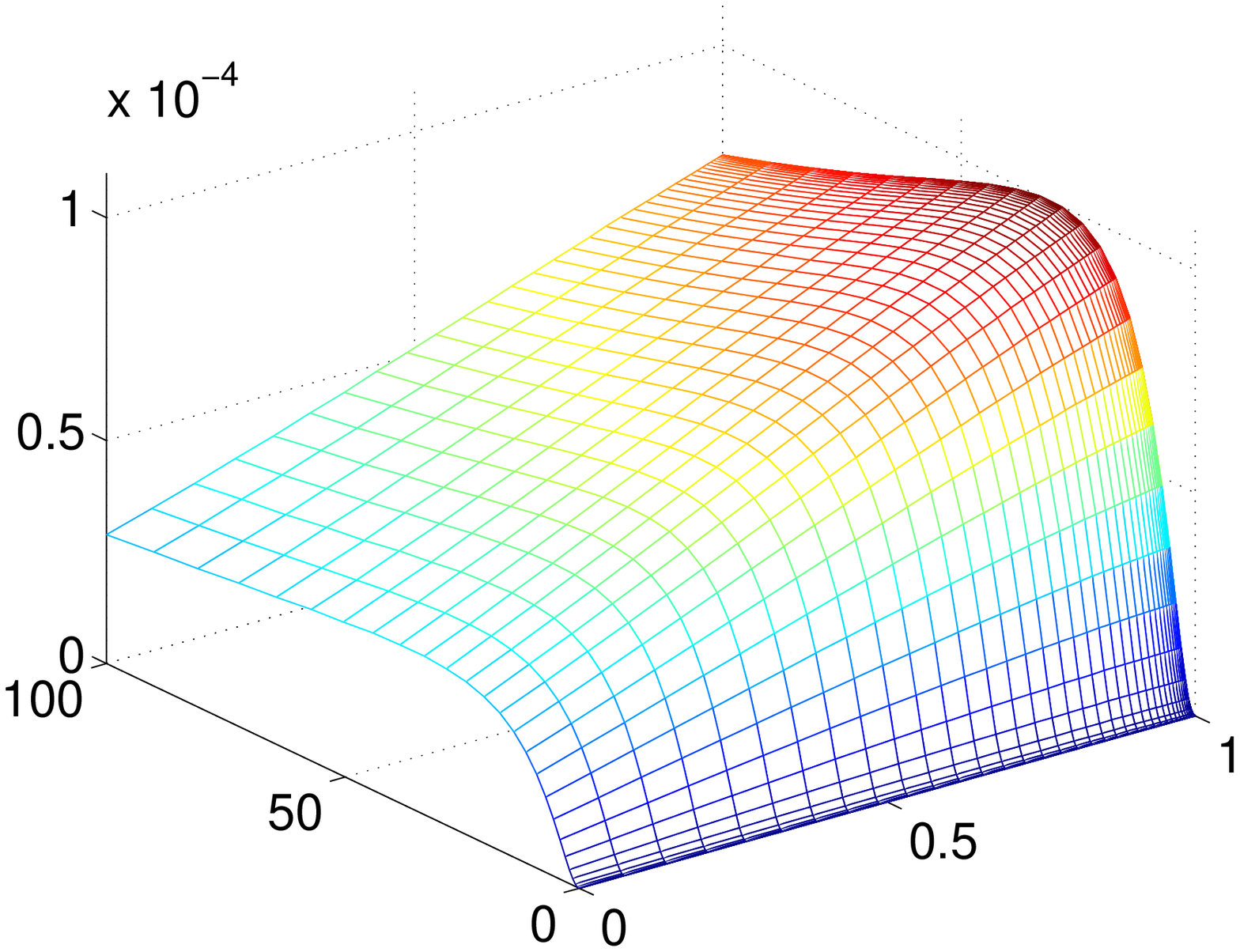}
    \put(-60,10){$x$}
    \put(-180,15){$t$}
    \put(-225,95){$\delta w$}
    \put(-225,150){b)}
    \caption{The relative error of the crack opening obtained for $N=40$ (nonuniform mesh, $\varrho=3$) and the time stepping strategy
    (\ref{t_step}) with $K=30$ .
    Fig.~\ref{blad_dw}a) corresponds to the \emph{solver 1}, while Fig.~\ref{blad_dw}b) refers to the \emph{solver 2}.}
\label{blad_dw}
\end{figure}

Let us now analyze the distribution of crack opening error, $\delta w$,  in time
and space. The respective results are presented in
Fig.~\ref{blad_dw}~a) -- \emph{solver 1}, and Fig.~\ref{blad_dw}~b)
-- \emph{solver 2}. In both cases the maximal errors are located at
the crack tip while the error distribution in time follows the trend
visible in Fig.~\ref{bl_czas}.

\begin{figure}[h!]
%M/N=1/300
    \hspace{6mm}\includegraphics [scale=0.40]{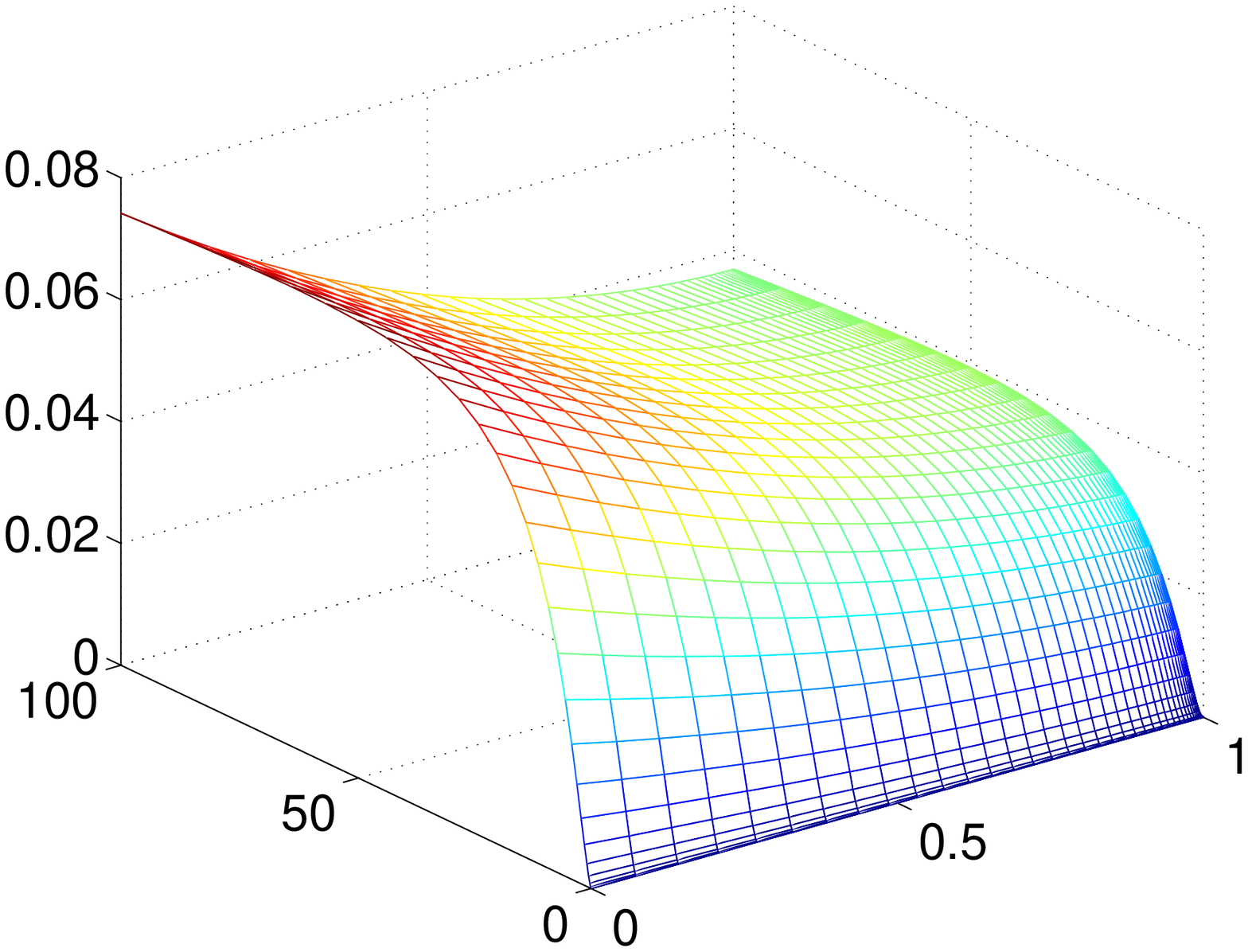}
    \put(-60,10){$x$}
    \put(-180,15){$t$}
    \put(-235,95){$\delta w_t$}
      \put(-225,150){a)}
    \hspace{-3mm}\includegraphics [scale=0.40]{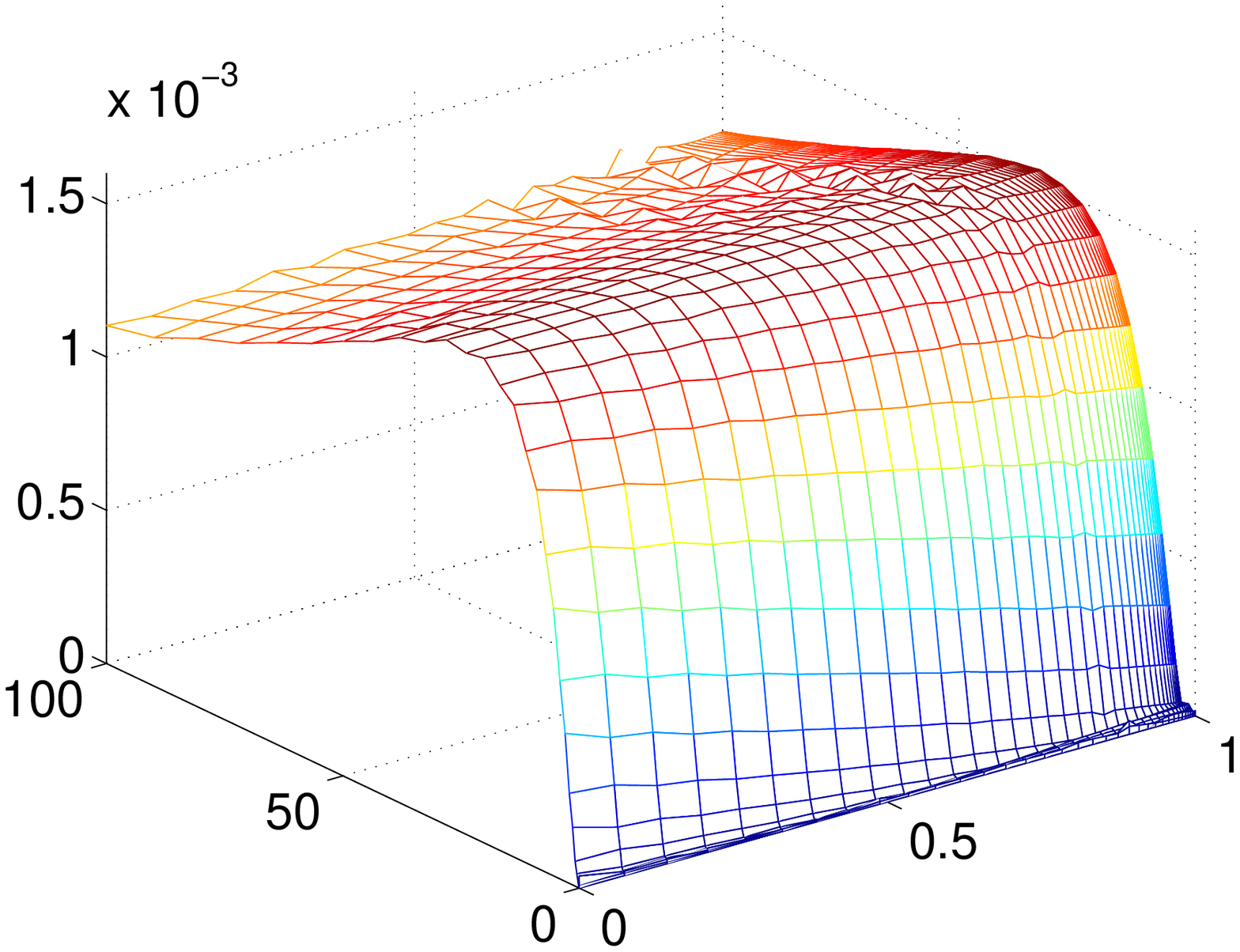}
    \put(-60,10){$x$}
    \put(-180,15){$t$}
    \put(-225,90){$\delta w_t$}
      \put(-225,150){b)}
    \caption{The relative error of the temporal derivative of the crack opening.
    Solution obtained by: a) \emph{solver 1}, b) \emph{solver 2} for $N=40$, non-uniform mesh ($\varrho=3$) and the time step strategy (\ref{t_step}) with $K=30$.}
\label{blad_dwt}
\end{figure}

Fig.~\ref{blad_dwt} shows the distributions of $\delta w_t$ obtained
by the integral solvers. It confirms our previous observation, that
\emph{solver 2} always provides better results than \emph{solver 1}.
Moreover, the greatest error in the case of \emph{solver 1} is located
at the crack inlet and the lowest at the crack tip, while
\emph{solver 2} gives approximately the same values of $\delta
w_t$ along the crack length.

\begin{figure}[h!]
%M/N=1/300
    \includegraphics [scale=0.40]{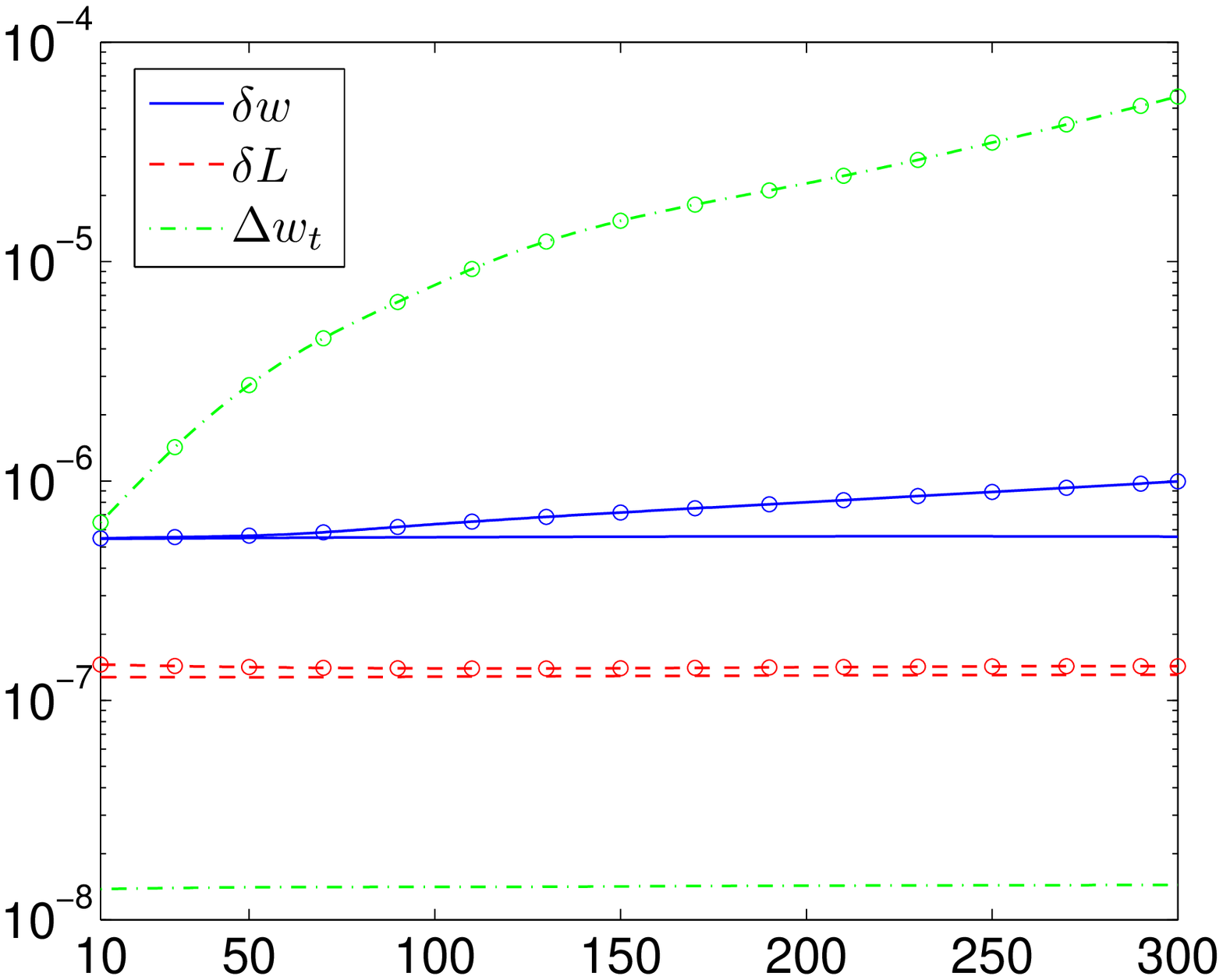}
    \put(-110,-2){$K$}
    \put(-225,150){a)}
    \hspace{-2mm}\includegraphics [scale=0.40]{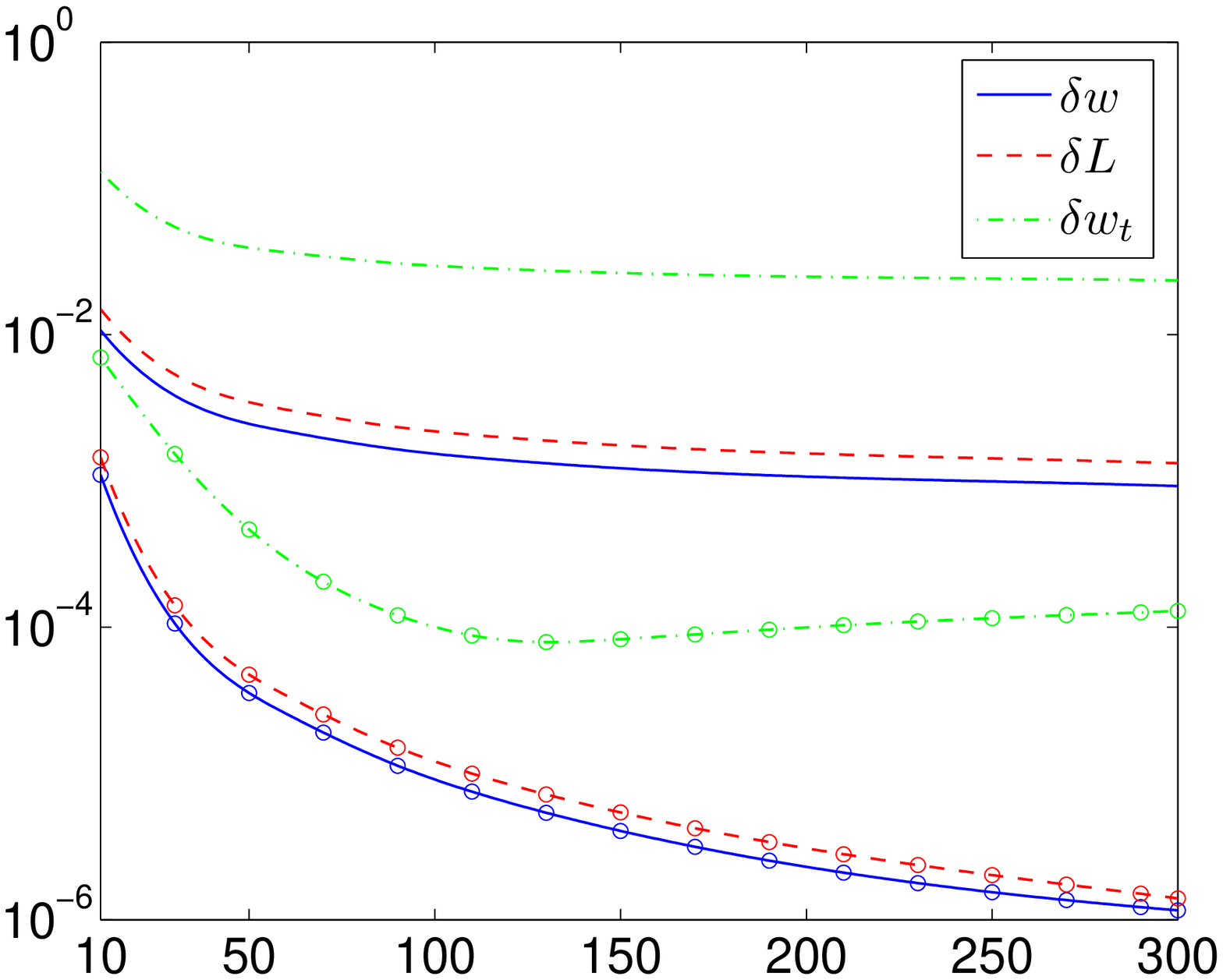}
    \put(-110,-2){$K$}
    \put(-225,150){b)}
    \caption{The errors of solution for two variants of $\gamma$:
    Fig.~\ref{blad_gamma}a) $\gamma=0$, Fig.~\ref{blad_gamma}b) $\gamma=1/3$. Line with markers correspond to \emph{solver 2}. }
\label{blad_gamma}
\end{figure}

In the last test in this subsection we  analyze the relation between the regimes of crack propagation and the performances of respective solvers. As mentioned previously, the benchmark solution in form \eqref{power} can be used to imitate various dynamic modes of the crack evolution. So far we have utilized the exponent of the time dependent term of the value $\gamma=1/5$ which refers to the constant injection flux rate. Now, let us consider two other variants of $\gamma$: i) $\gamma=0$ - for this choice the normalized crack opening is constant in time; ii) $\gamma=1/3$ - this value corresponds to the steady state propagation of the fracture. For the computations the same spatial mesh as before is taken ($N=40$). The number of time steps accepted within strategy \eqref{t_step}, $K$, ranges from 10 to 300. In this way the graphs (Fig.~\ref{blad_gamma}a)-b)) describing the solution errors in the function of $K$ were prepared for both values of $\gamma$ (similarly as in Fig.~\ref{bl_dt}). The analyzed accuracy parameters were: relative error of the crack opening, $\delta w$, relative error of crack the length, $\delta L$, relative error of the crack opening temporal derivative, $\delta w_t$, for $\gamma=1/3$ and absolute  error of the crack opening temporal derivative, $\Delta w_t$, for $\gamma=0$.

The results depicted in Fig.~\ref{blad_gamma} show that for $\gamma=0$ one obtains much more accurate results than for $\gamma=1/3$ that could have been predicted ($w_t=0$). However, it is a surprise that for $\gamma=0$ \emph{solver 1} provides better solution accuracy than \emph{solver 2}. Although the difference is moderate in case of $\delta w$ and $\delta L$, the values of $\Delta w_t$ vary by at least two orders of magnitude. From Fig.~\ref{blad_gamma}a) it follows that for this regime of crack propagation the solution accuracy cannot be improved by simple refining the temporal mesh, and for \emph{solver 2} even a reverse relation is observed.

The situation is quite different for $\gamma=1/3$ (Fig.~\ref{blad_gamma}b) ). This time again \emph{solver 2} proves its advantage over \emph{solver 1} for all the analyzed parameters. For \emph{solver 2}, it is sufficient to take only 30 time steps to have much better results than those provided by \emph{solver 1} for 300 steps. The solution accuracy can be improved by increasing the number of time steps, however it seems that for \emph{solver 1} the saturation level is close to $K=300$. A similar trend was observed for $\gamma=1/5$ (see Fig.~\ref{bl_dt}).

A direct conclusion from this test is that for different modes of crack propagation there are different optimal time stepping strategies. This should properly accounted for especially in the cases when the values of injection flux rate or leak-off to formation change appreciably in the considered time interval.

The aforementioned analysis prove that in terms of accuracy in most of the cases
\emph{solver 2} is much better that \emph{solver 1} with respect to
all computed components of the solution: the crack length, $L$, the
crack opening, $w$, its temporal derivative, $w_t$ and the fracture propagation speed, $V_0$. However, for some regimes of crack propagation (low values of $\gamma$)  \emph{solver 1} may give comparable or even slightly better results than  \emph{solver 2}.
The advantage of \emph{solver 1} is better efficiency of
computations: the time of computations for this solver was on
average one third lower that for \emph{solver 2}.

\subsubsection{Example with singular leak-off regime.}
\label{Carter}

In \eqref{q_l} we have assumed that the behaviour of the leak-off
function near the crack tip can be described by a power law, giving
in the worst case a square root singularity. Such a limiting
behaviour corresponds to the Carter leak-off model \cite{Carter}. As
a result, although the leading term of the asymptotic expansion for
the crack opening near the fracture tip remains the same, the higher
terms change, disturbing the solution smoothness (see
\cite{Kovalyshen}). A comprehensive analysis of this case was done
in \cite{Kusmierczyk} where it was also proved that the
deterioration of solution accuracy can be prevented by employing the
second asymptotic term in the computational algorithm.

In this subsection we show that the algorithms developed in the
paper are capable of tackling this kind of problems \emph{without
any additional modifications}. To this end let us consider another
benchmark solution $u(t,x)$ (see (\ref{u_bench})) defined by the
functions:
\begin{equation}
\label{Carter_bench}
h(x)=(1-x)^{1/3}(1+s(x)),\quad
s(x)=\frac{1}{5}(1-x)^{1/6},
\end{equation}
with the same function $\psi(t)$ with $\gamma=1/5$.

One can easily check, that the above form of $s(x)$ results in a
singular behaviour of $q_l$, with the leading term  of the order
$O\big((1-x)^{-1/2}\big)$ as $x\to1$. The value of multiplier $u_0$
in \eqref{u_bench} was taken in such a way to make the benchmark
comparable with the one used previously, in a sense of an average
particle velocity. Indeed, in \cite{M_W_L} for a fluid velocity
defined as $V=q/w$, a parameter describing its variation along the
crack length was introduced:
\begin{equation}
\label{gamma_v}
\gamma_v(t)=\big[\max_{x}\left(V(x,t)\right)-\min_{x}\left(V(x,t)\right)\big]\left[\int_0^1V(x,t)dx\right]^{-1}.
\end{equation}
This parameter reflects indirectly the balance between the flux
injection rate and leak-off to formation. It was also shown there
that it has a decisive influence on the accuracy of computations
(the greater value of $\gamma_v$, the greater error of the
computations). For the benchmarks with comparable values of
$\gamma_v$, one can expect similar accuracy of the computations.
This trend was also confirmed in \cite{Kusmierczyk}. In our case,
the deterioration of the solution smoothness near the crack tip is an
additional factor which contributes to the increase of the computational
error.

Note that because of the chosen structure of (\ref{u_bench}), the
value of $\gamma_v$ is constant in time for all benchmarks
considered in this paper. For the benchmark \eqref{s_1} $\gamma_v$
yields 0.408, while \eqref{Carter_bench} one leads to
$\gamma_v=0.411$.

\begin{figure}[h!]
%M/N=1/300
    \hspace{-2mm}\includegraphics [scale=0.40]{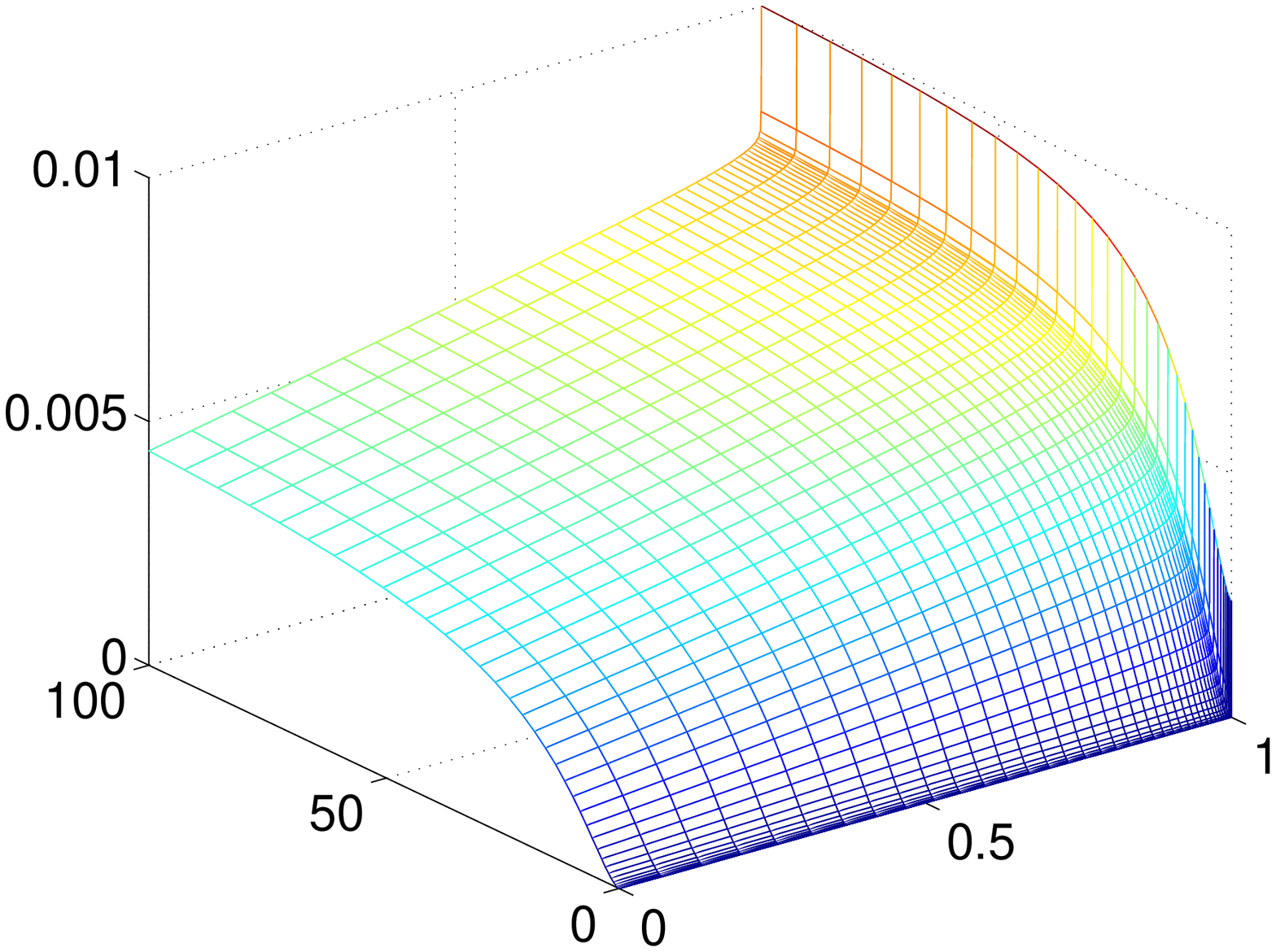}
    \put(-60,10){$x$}
    \put(-180,15){$t$}
    \put(-220,110){$\delta w$}
    \hspace{-6mm}\includegraphics [scale=0.40]{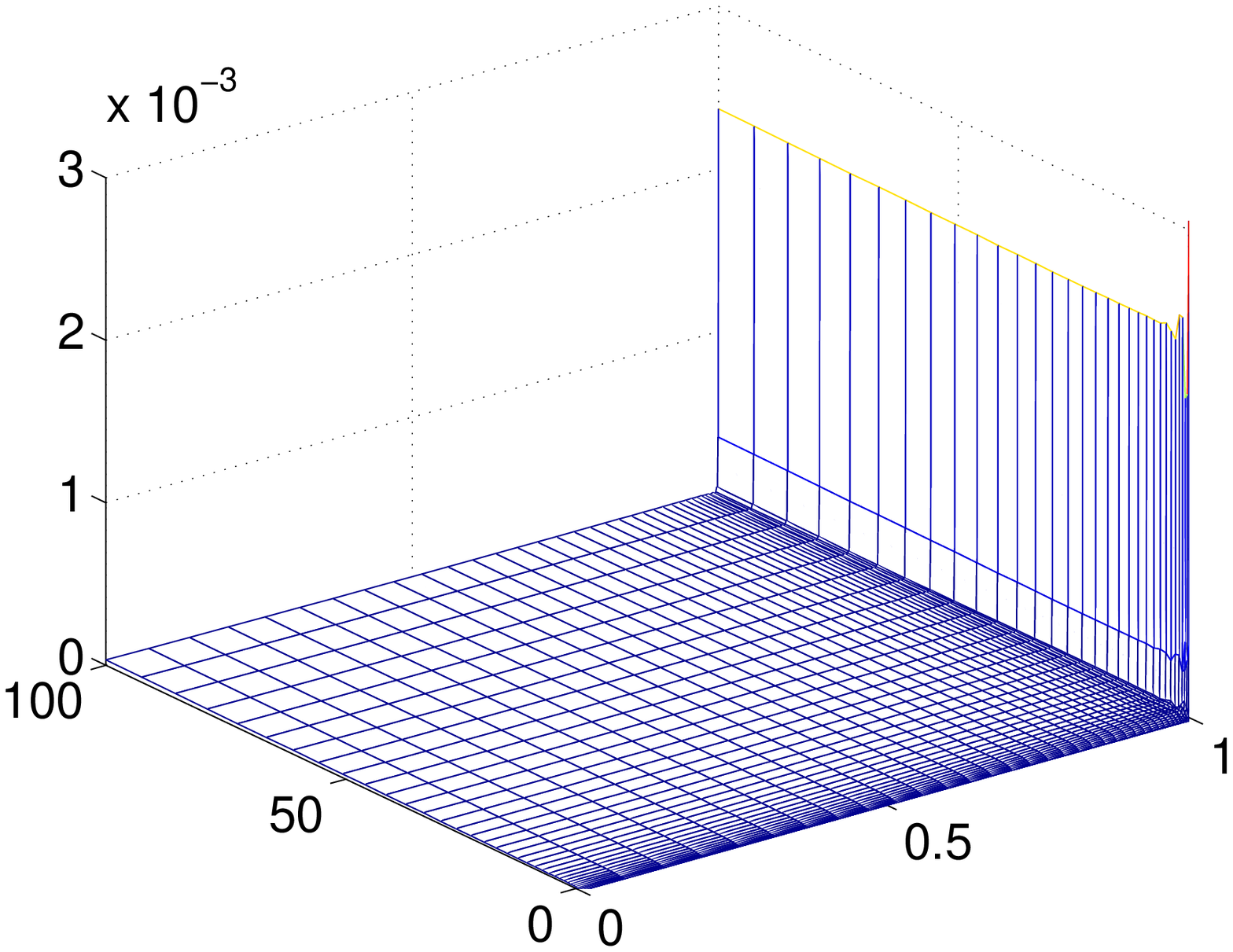}
    \put(-60,10){$x$}
    \put(-180,15){$t$}
    \put(-220,95){$\delta w$}
    \caption{The relative error of the crack opening. Solution obtained by: a) \emph{solver 1}, b) \emph{solver 2}. Other parameters in the computations were: $N=40$, $K=30$, $\varrho=3$.}
\label{blad_dw_carter}
\end{figure}

\begin{figure}[h!]
%M/N=1/300
    \hspace{4mm}
    \includegraphics [scale=0.40]{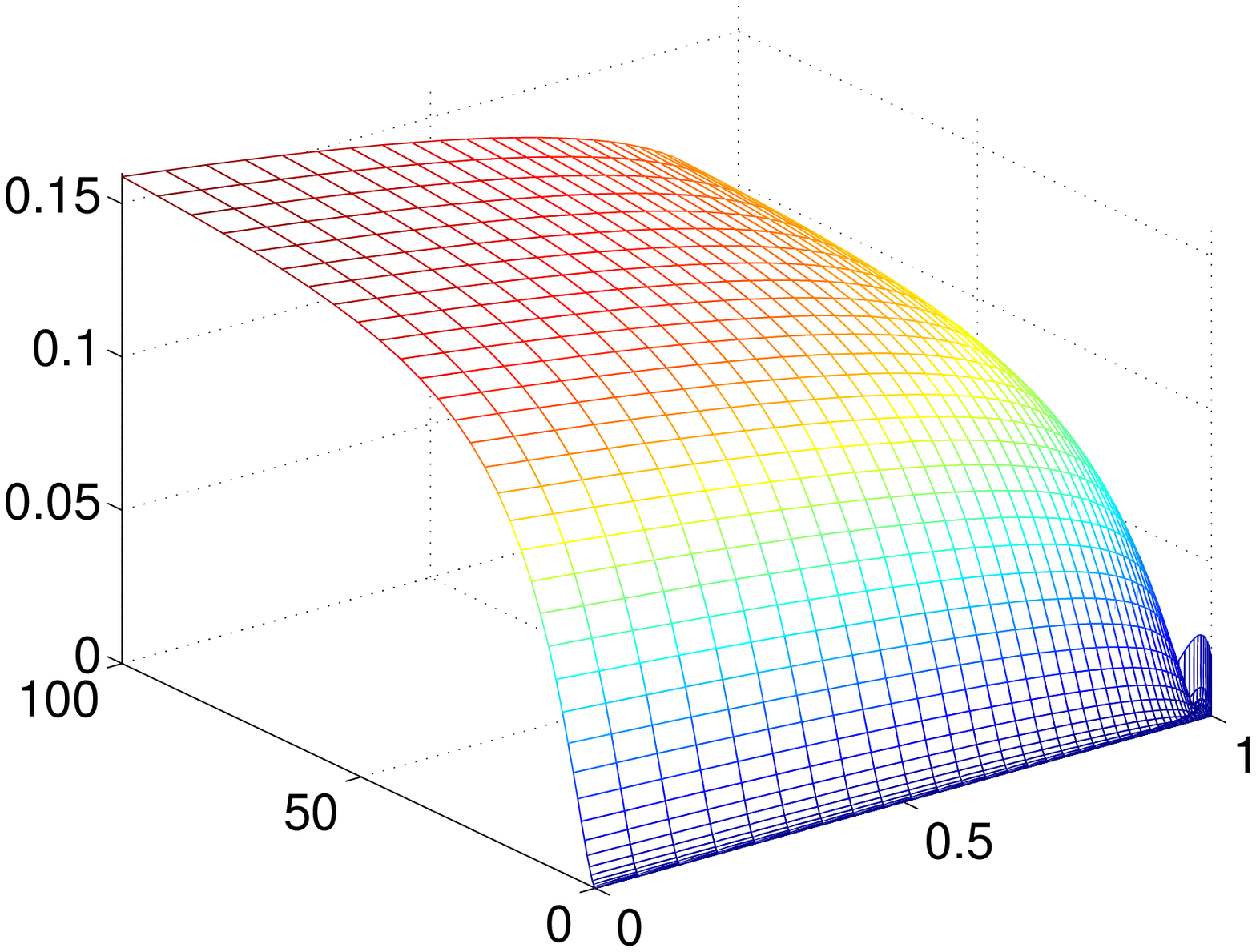}
    \put(-60,10){$x$}
    \put(-180,15){$t$}
    \put(-235,95){$\delta w_t$}
    \put(-230,150){a)}
    \hspace{2mm}\includegraphics [scale=0.40]{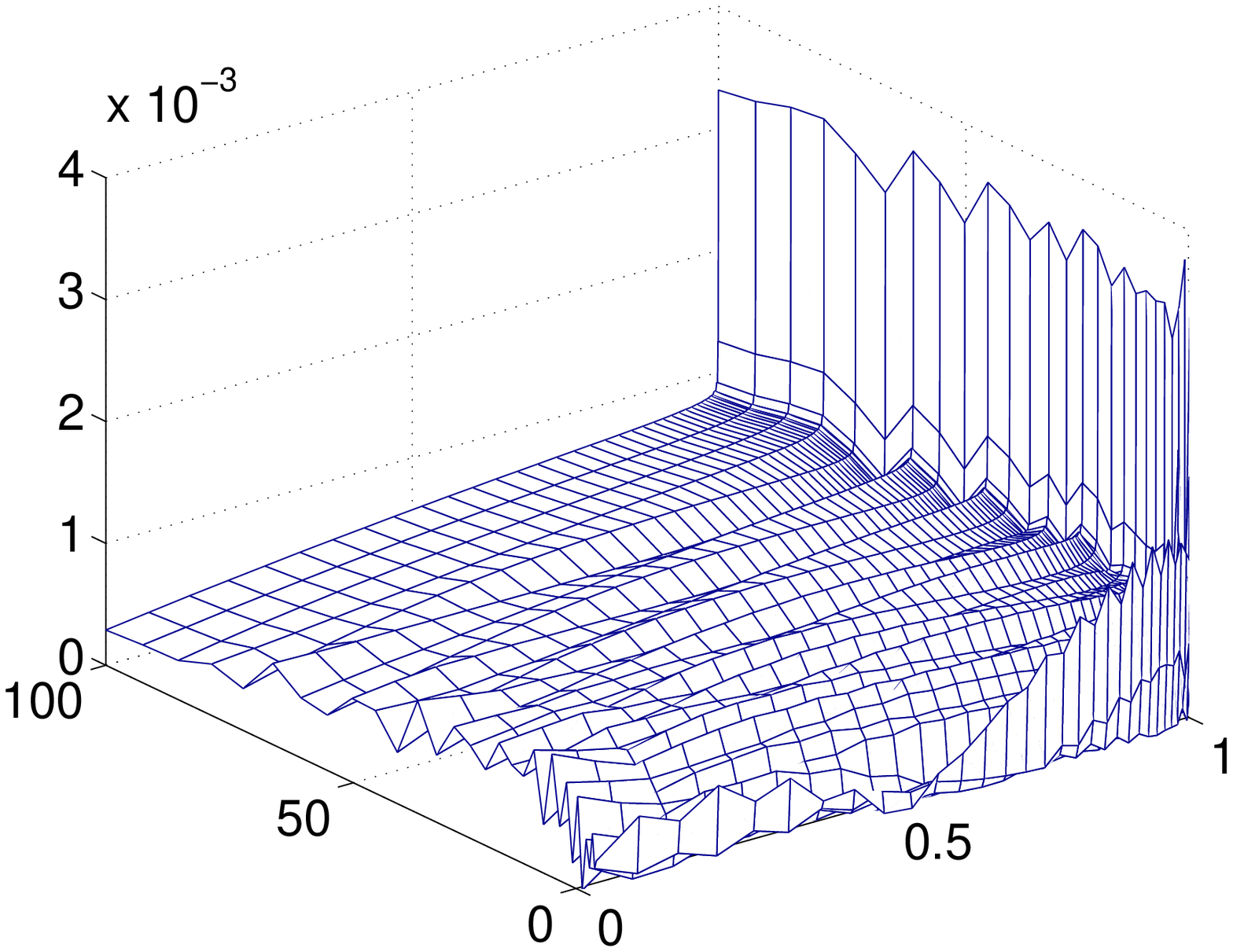}
    \put(-60,10){$x$}
    \put(-180,15){$t$}
    \put(-225,90){$\delta w_t$}
    \put(-230,150){b)}
    \caption{The relative error of the temporal derivative of crack opening. Solution obtained by: a) \emph{solver 1}, b) \emph{solver 2}.
    Other parameters in the computations are: $N=40$, $K=30$, $\varrho=3$.}
\label{blad_dwt_carter}
\end{figure}

The computations were done for the same spatial
($N=40$, $\varrho=3$) and temporal ($K=30$, strategy \eqref{t_step})
meshes as previously considered. The distributions of the errors $\delta w$ and $\delta w_t$ are shown in
 Fig.~\ref{blad_dw_carter} -- Fig.~\ref{blad_dwt_carter}, respectively.
Comparing these results with those obtained for the finite leak-off
regime (see Fig.~\ref{blad_dw} -- Fig.~\ref{blad_dwt}) one can see
that the accuracy of computations decreased significantly.

In the case of \emph{solver 1}, the crack opening error, $\delta w$,
shown in Fig.~\ref{blad_dw_carter} a) is of one order of magnitude
greater than that reported in Fig.~\ref{blad_dw} a). For
\emph{solver 2} in turn, $\delta w$ depicted in
Fig.~\ref{blad_dw_carter} b) is almost two orders higher than the
one obtained previously for $q_l$ vanishing at $x=1$
(Fig.~\ref{blad_dw} b)). A pronounced jump of the error is observed
at the crack tip, especially for \emph{solver 2}, which explains the
deterioration of accuracy when comparing with the finite leak-off
case. However, if one considers the accuracy of the solution away from
the crack tip, it is still of the same quality as before.

The same trend can be observed for the temporal derivative of the
crack opening, $w_t$, see Fig.~\ref{blad_dwt_carter}. Opposite to
the benchmark case \eqref{s_1} with the finite leak-off, the
distribution of $\delta w_t$ for \emph{solver 2} becomes highly
non-uniform, with distinct increase at the crack tip.

This deterioration of the solution accuracy near $x=1$ for the Carter
leak-off model should not be a surprise, as the algorithm described
above in (\ref{dm1}) --  (\ref{G_2}) accounted directly only  for
the first (leading) term of the asymptotic expansion for the crack
opening. For the singular leak-off model, the function $\Delta u$
and other integrands employed in  $G_1$ and $G_2$ are not
sufficiently smooth near the crack tip, which increases the errors
of integration. However, one can counteract this tendency by
accounting for further asymptotic terms in the algorithm, which will
be illustrated in the following with the example of two terms of $w$
expansion.

\begin{figure}[h!]
%M/N=1/300
    \hspace{6mm}
    \\[2mm]\includegraphics [scale=0.40]{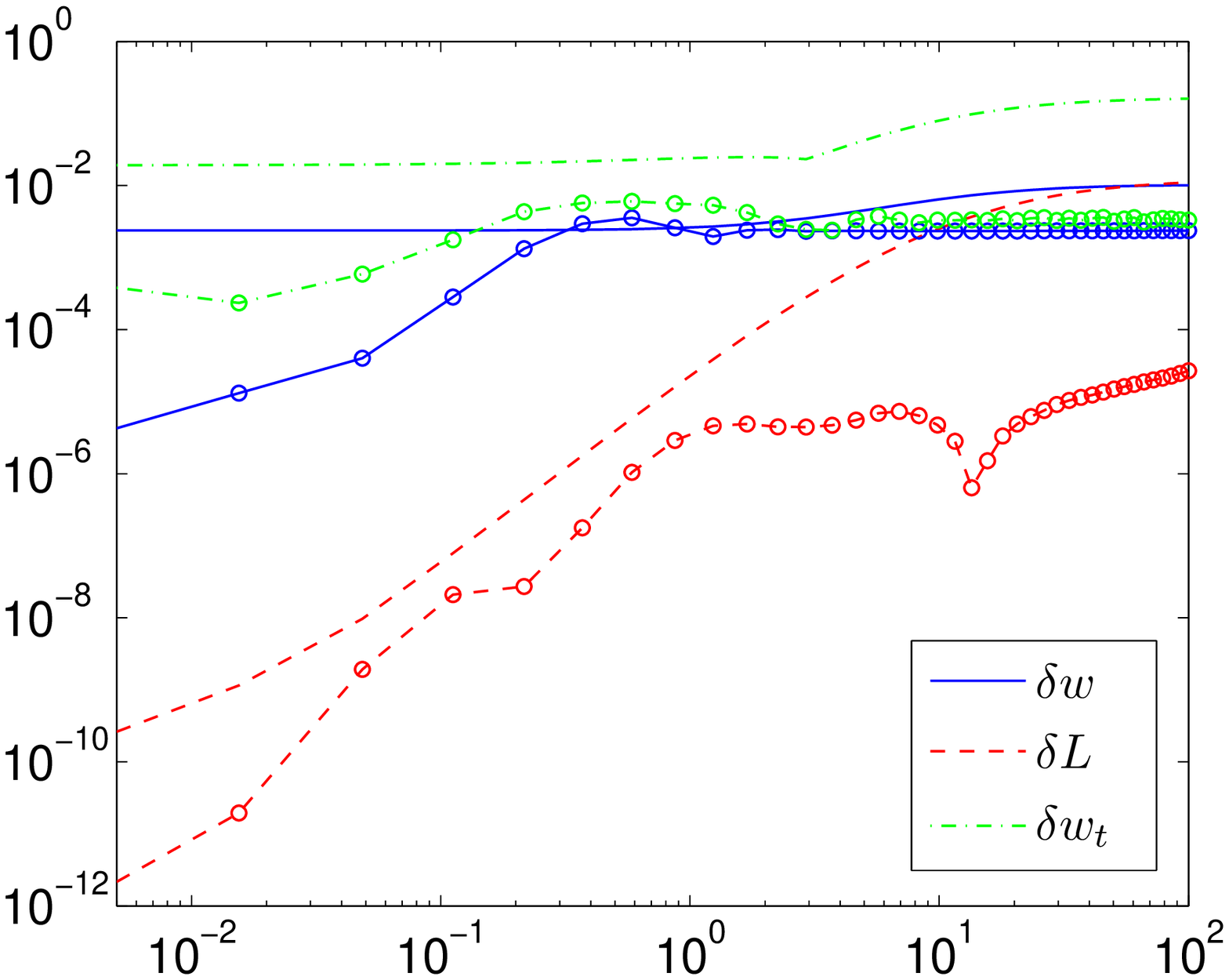}
    \put(-110,0){$t$}
    \put(-230,150){a)}
    \hspace{-2mm}
    \includegraphics [scale=0.40]{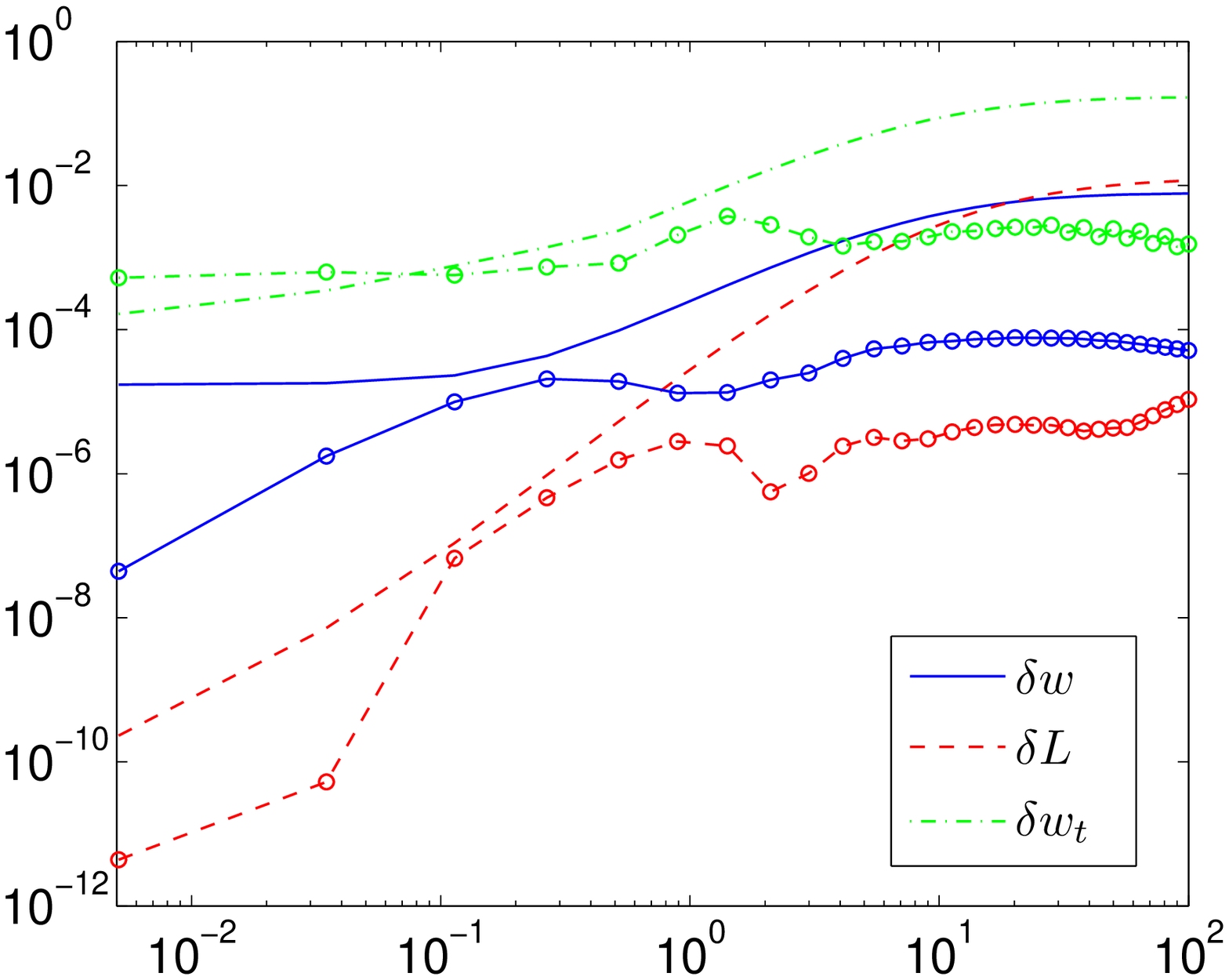}
    \put(-110,0){$t$}
    \put(-230,150){b)}
    \caption{The relative error of the temporal derivative of crack opening. Solution obtained by: a) original algorithms, b) modified algorithms.
    The lines with markers refer to \emph{solver 2}}
\label{blady_carter}
\end{figure}

Let us investigate the evolution of solution errors ($\delta w$,
$\delta L$ and $\delta w_t$) in time in the way it was done in
Fig.~\ref{bl_czas} for the case of the finite leak-off. In
Fig.~\ref{blady_carter} a) we show the results obtained by direct
execution of the algorithm (\ref{dm1}) -- (\ref{G_2}).
Fig.~\ref{blady_carter} b) depicts the case of a modified algorithm
employing also the second asymptotic term of the crack opening. When
analyzing the results, one observes again the trend of $\delta L$
accumulation for both solvers. Interestingly, $\delta L$ has the
same level for both the original and the modified algorithms. In
the case of \emph{solver 1} for other analyzed parameters as well only a
very slight improvement is reported. However, for \emph{solver 2}
one obtains a pronounced reduction of computational errors for $w$
and $w_t$. To have a further advance here, the next terms of
asymptotic expansion of the crack opening may be taken into account.
On the other hand, for the Carter leak-off model the number of nodal
points for which the accuracy saturation is observed is larger than
for a finite leak-off. Thus, if necessary, one can also improve the
accuracy  by increasing the value of $N$ and by taking a greater
number of time steps $K$.

From the analysis given in this subsection two main conclusions can
be drawn. First, that the problems with singular leak-off regimes
can be directly tackled by the proposed algorithms. Also in such
cases \emph{solver 2} yields more accurate results than \emph{solver
1}. The second conclusion is that in order to have the solution
accuracy comparable to that achieved for non-singular leak-off, one
has to use a larger number of nodal points and/or employ further terms of asymptotic expansion for $w$ in the
algorithm.

\section{Conclusions}

We would like to itemize the following conclusions as a resume of
this paper:
\begin{itemize}
\item Presented approach can be efficiently used for tackling the PKN model of hydrofracturing and may be adopted for multifracture systems.
\item Both new solvers provide better
computational accuracy  than the conventional algorithms from \cite{M_W_L}. Moreover, comparable accuracies can be achieved
here at much lower computational cost, as the new solvers enable us to drastically reduce the
densities of spatial and temporal meshes.

\item New solvers are appropriate for directly tackling the problems with different fluid flow regimes, including various injection flux rates and singular leak-off.
\item
In order to increase the efficiency and accuracy of computations, it is advisable to employ at least two asymptotic terms of the crack opening, $w$.
\item The developed algorithms do not require any regularization techniques.
The boundary conditions are imposed directly into the numerical scheme. The speed equation plays a crucial
role in the analysis.
\end{itemize}

\subsubsection*{Acknowledgements}
This work has been done in the framework of the EU FP7 PEOPLE project under contract
number PIAP-GA-2009-251475-HYDROFRAC. The authors are grateful to the
Institute of Mathematics and Physics of Aberystwyth University and
EUROTECH Sp. z o. o. for the facilities and hospitality.

\section*{Appendix A: Numerical benchmarks}

Let us define a set of benchmark solutions useful for testing
different numerical solvers. Consider a class of positive functions
$C_+(0,1)$ described in the following manner:
\[
C_+(0,1)=\{h\in C^2(0,1)\cap C[0,1],
\]
\[  \lim_{x\to 1-}
(1-x)^{-1/3}h(x)=1,\quad h(x)> 0, \,\,\, x\in[0,1)\}.
\]
By taking an arbitrary $h\in C_+(0,1)$, one can build a benchmark
solution for the normalized formulation of the problem as:
\begin{equation}
\label{u_bench} u(x)=u_0\psi_j(t)h(x).
\end{equation}
where functions $\psi_j(t)$ and $h(x)$ are specified below. On
substitution of \eqref{u_bench}  into  (\ref{Rey1}) one finds:
\begin{equation}
\label{q_l_bench} q_l(t,x)=
\end{equation}
\[
\gamma u_0\big[\frac{1}{\beta}\left(x
h^3(x)h'(x)+3\left(h^3( x)h'(x)\right)'\right)- h(
x)\big]\psi_j^\alpha(t),
\]
\[
w_0(t)=u_0\psi_j(t),
\]
where  two sets of the benchmark solutions can be considered. For
the first one we choose $\psi_1(t)=e^{\gamma t}$ and $\beta=2/3$,
$\alpha=1$, while for the second, $\psi_2(t)=(a+t)^\gamma$ and
$\beta=2\gamma/(3\gamma+1)$, $\alpha=(3\gamma-1)/\gamma$.
Corresponding crack lengths are defined in (\ref{L_exp}) and
(\ref{L_power}), respectively. Finally, the injection flux rate is
computed from the boundary condition \eqref{BC_11}$_1$:
\begin{equation}
\label{q_0_bench} q_0(t)=-\frac{u_0^4\psi^4_j(t)}{L(t)}h^3(0)h'(0),
\end{equation}
while the initial condition reads $W(x)=u_0\psi_j(0)h(x)$.

Note that when taking the function $h(x)$ from the class $C_+(0,1)$
in the following representation:
\[
h(x)=(1-x)^{1/3}(1+s(x)),
\]
\[s\in C^2[0,1], \quad s( x)>-1,\,\,\,
x\in[0,1)
\]
$q_l$ automatically satisfies the condition $
q_l(t,x)=O\big((1-x)^{1/3}\big),\quad x\to1. $

The presented benchmarks allow one to test numerical schemes in various
fracture propagation regimes (accelerating/decelerating cracks) by
choosing proper values of the parameter $\gamma$ (see \cite{M_W_L}).
Additionally, if one reduces the requirements for the smoothness of
 $s(x)$ near $x=1$, assuming that
$h\in C^2[0,1)\bigcap H^\alpha (0,1)$ the benchmark can serve to
model singular leak-off regimes (compare with (\ref{Carter_bench})).

Note that the zero leak-off case cannot be described by the
aforementioned group of benchmarks. However, an analytical benchmark
for this regime, represented in terms of a rapidly converging series,
has been developed in \cite{Linkov_3}.


\begin{thebibliography}{99}



\bibitem{Adachi_2002} Adachi J, Detournay E (2002) Self-similar solution of a plane-strain fracture driven by a power-law fluid.
Int J Numer Anal Methods Geomech 26: 579-604

\bibitem{Adachi_2008} Adachi J, Detournay E (2008) Plane strain propagation of a hydraulic fracture in a permeable rock. Eng Fract Mech 75(16): 4666-4694

\bibitem{Adachi_2007b} Adachi J, Peirce A (2007) Asymptotic analysis of an elasticity equation for a finger-like hydraulic fracture. J Elast 90(1): 43-69

\bibitem{Adachi_2007} Adachi J, Siebrits E, Peirce A, Desroches J
(2007) Computer simulation of hydraulic fractures.  Int J
Rock Mech Min Sci 44: 739-757


\bibitem{Carter} Carter E (1957) Optimum fluid characteristics for fracture extension. In: Howard, G., Fast, C. (eds.) Drilling and
Production Practices, 261-270. American Petroleum Institute

\bibitem{Crittendon} Crittendon BC (1959) The mechanics of design and
interpertation of hydraulic fracture treatments. J Pet
Tech 21: 21-29


\bibitem{Detournay}Detournay E (2004) Propagation regimes of fluid-driven fractures in impermeable rocks. Int J Geom 4:
1-11

\bibitem{Dobroskok} Dobroskok AA, Linkov AM (2011) Modeling of fluid flow, stress state and seismicity induced in rock by instant pressure drop in a hydrofracture. J Min Sci 47(1): 10-19

\bibitem{Economides2000}  Economides M, Nolte K (eds.) (2000) Reservoir
Stimulation. 3rd edn. Wiley, Chichester, UK

\bibitem{Garagash et al}
 Garagash DI,  Detournay E (2000)
The tip region of a fluid-driven fracture in an elastic
medium. J Appl Mech 67: 183-192

\bibitem{Geertsma} Geertsma J, de Klerk F (1969) A rapid method of
predicting width and extent of hydraulically induced fractures.
J Pet Tech 21: 1571-1581 [SPE 2458]


\bibitem{Hubbert} Hubbert MK, Willis DG (1957) Mechanics of hydraulic
fracturing. J. Pet. Tech. 9(6): 153-68

\bibitem{Kemp} Kemp LF (1989) Study of Nordgren's equation of hydraulic fracturing. SPE Production Eng 5:
311-314

\bibitem{Khristianovic} Khristianovic SA, Zheltov YP (1955) Formation of
vertical fractures by means of highly viscous liquid. In:
Proceedings of the fourth world petroleum congress, Rome, 1955,
579-586

\bibitem{Kovalyshen} Kovalyshen Y, Detournay E (2009) A reexamination
of the classical PKN model of hydraulic fracture.  Transp Porous
Med 81: 317-339

\bibitem{Kovalyshen_1} Kovalyshen Y (2010) Fluid-driven fracture in poroelastic
medium. PhD Thesis, The University of Minnesota


\bibitem{Kusmierczyk} Kusmierczyk P, Mishuris G, Wrobel M (2012) Remarks on numerical simulation of the PKN model of hydrofracturing in proper variables. Various leak-off regimes. arXiv:1211.6474.

 \bibitem{Lecampion_Det}
 Lecampion B,  Detournay E  (2007) An implicit algorithm for the propagation of a hydraulic fracture with a
fluid lag. Comput Method Appl M 196(49-52): 4863 -- 4880

\bibitem{Linkov_1} Linkov AM (2011) Speed equation and its application
for solving ill-posed problems of hydraulic fracturing. {\it ISSM
1028-3358, Doklady Physics}, 56(8): 436-438. Pleiades Publishing,
Ltd.

\bibitem{Linkov_2} Linkov AM (2011) Use of a speed equation for
numerical simulation of hydraulic fractures. arXiv:1108.6146

\bibitem{Linkov_3} Linkov AM (2011) On efficient simulation of
hydraulic fracturing in terms of particle velocity. Int J
Engng Sci 52: 77-88

\bibitem{M_W_L} Mishuris G, Wrobel M, Linkov A (2012) On modeling hydraulic fracture in proper variables:
stiffness, accuracy, sensitivity. Int J Engng Sci
61: 10-23

\bibitem{Mitchell_2007} Mitchell SL, Kuske R, Peirce A (2007) An asymptotic framework for finite hydraulic fractures including leakoff.
SIAM J Appl Math 67(2): 364--386

\bibitem{Mitchell_2007b} Mitchell SL, Kuske R, Peirce AP (2007) An asymptotic framework for the analysis of hydraulic fractures: the impermeable case. ASME J Appl Mech  74(2): 365--372


\bibitem{Moos} Moos D (2012) The importance of stress and fractures in hydrofracturing and stimulation performance: a geomechanics overview. Search and Discovery Article 80255 (2012)

\bibitem{Nordgren} Nordgren, RP (1972) Propagation of a Vertical
Hydraulic Fracture.  J Pet Tech 253: 306-314


\bibitem{Olson} Olson JE (2008) Multi-fracture propagation modeling: Applications to hydraulic fracturing in shales and tight gas sands. The 42nd U.S. Rock Mechanics Symposium (USRMS), June 29 - July 2, 2008, San Francisco, CA

\bibitem{Pierc_Det} Peirce A, Detournay E (2008). An implicit level set method for modeling hydraulically driven fractures. Comput Methods Appl Mech Engrg  197: 2858--2885
\bibitem{Pierc_Det_2}
    Peirce A, Detournay E (2009), An eulerian moving front algorithm with weak-form tip asymptotics for
modeling hydraulically driven fractures. Num Meth Eng 25(2): 185-200

\bibitem{Perkins} Perkins TK, Kern LR (1961) Widths of hydraulic fractures. J
Pet Tech  13(9): 37-49 [SPE 89]

\bibitem{Rubin} Rubin AM (1995) Propagation of magma filled cracks. Ann
Rev Earth Planet Sci 23: 287-336

\bibitem{Sneddon} Sneddon IN, Elliot HA (1946) The opening of a Griffith
crack under internal pressure. Q Appl Math 4:
262-267

\bibitem{stiff_system} Stoer J, Bulirsh R (2002) Introduction to
Numerical Analysis. Third Edition. Springer-Verlag.

\bibitem{Rice} Tsai VC, Rice JR (2010) A model for turbulent
hydraulic  fracture and application to crack propagation at
glacier beds. J Geophys Res 115: 1-18


\bibitem{Zhang} Zhang X, Jeffrey R, Llanos EM (2004) A study of shear hydraulic fracture propagation. Gulf Rocks 2004, the 6th North America Rock Mechanics Symposium (NARMS), June 5 - 9, 2004 , Houston, Texas


\end{thebibliography}
\end{document}